\documentclass{article}
\usepackage{fullpage}
\usepackage{amsthm}
\usepackage{amsfonts}
\usepackage{amsmath}
\usepackage{rotating}
\usepackage{gnuplot-lua-tikz}
\usepackage{url}

\numberwithin{equation}{section}

\usepackage{tikz}
\usetikzlibrary{matrix,calc,arrows,decorations.markings}
\tikzstyle{vecArrow} = [
   line width=4,
   decoration={markings,mark=at position 1 with {\arrow[thick]{triangle 60}}},
   shorten >= 5.5pt,
   preaction = {decorate},
]

\usepackage{listings}
\lstdefinelanguage{pseudocode}
{
   morekeywords={do,while,if,else,then,end,cycle,function,return},
   morecomment=[s]{/*}{*/},
   morecomment=[l]{//}
}
\newcommand{\pmat}[1]{\begin{pmatrix}#1\end{pmatrix}}

\newcommand\T{\rule{0pt}{2.6ex}}       % Top strut in table
\newtheorem{lemma}{Lemma}[section]
\newtheorem{theorem}{Theorem}[section]

\begin{document}
   \title{
      Compressed threshold pivoting for sparse symmetric indefinite systems
   }
   \author{
      J.\ D.\ Hogg and J.\ A.\ Scott
   }
   \maketitle

   \begin{abstract}
      A key technique for controlling numerical stability in sparse direct
      solvers is threshold partial pivoting. When selecting a pivot, the entire
      candidate pivot column below the diagonal must be up-to-date and must be
      scanned. If the factorization is parallelized across a large number of
      cores, communication latencies can be the dominant computational cost.

      In this paper, we propose two alternative pivoting strategies for sparse
      symmetric indefinite matrices that
      significantly reduce communication by compressing the necessary data into
      a small matrix that can be used to select pivots. Once pivots have been
      chosen, they can be applied in a communication-efficient fashion.

      For an $n\times p$ submatrix on $P$ processors, we show our methods
      perform a factorization using $O(\log P)$ messages instead of the
      $O(p\log P)$ for threshold partial pivoting. The additional
      costs in terms of operations and communication bandwidth are relatively
      small.

      A stability proof is given and numerical results using a range of
      symmetric indefinite matrices arising from practical problems are used to
      demonstrate the practical robustness. Timing results on large random
      examples illustrate the potential speedup on current multicore
      machines.
   \end{abstract}

   \section{Introduction}

   We are interested in the efficient and stable factorization of large sparse
   symmetric indefinite matrices.
   Most algorithms for this employ supernodes (see, for example,
   \cite{adek:2001,gujk:2001,hosc:2010c,hosc:2011a,Reid2008,scga:2006}).
   That is, a set of consecutive
   columns having the same (or similar) sparsity pattern in the factor.
   By storing only those rows that contain nonzeros, each supernode
   may be held as a dense $n\times p$ trapezoidal matrix $A$ of the form
   \medskip
   \begin{center}
      \begin{tikzpicture}
         \draw (-0.5,0) -- (0.5,0) -- (-0.5,1) -- cycle; 
         \draw (0.5,0) -- (0.5,-4) -- (-0.5,-4) -- (-0.5,0);
         \node[inner sep=0] (a0) at (0.7,0.5) {};
         \node at (-0.1,0.3) {$A_{11}$};
         \node at (0,-2.0) {$A_{21}$};
         \node (aA21) at (0.5,-2.0) {};
         \draw[<->] (-0.7,1) -- (-0.7,-4);
         \draw[<->] (-0.5,1.1) -- (0.5,1.1);
         \node[fill=white,font=\footnotesize,outer sep=0pt,inner sep=1pt] at (-0.7,-1.5) {$n$};
         \node[fill=white,font=\footnotesize,outer sep=0pt,inner sep=1pt] at (0.0,1.1) {$p$};
         \node at (-1.3,-1.5) {$A=$};
 
      \end{tikzpicture}
      \vspace*{-3.00cm}\begin{eqnarray}\label{supernode}\end{eqnarray}\vspace*{0.0cm}
   \end{center}
   \vspace{1.5cm}
This matrix is termed a {\it supernodal} matrix. In general, at non-root nodes, $n\gg p$.
   \newpage
   \noindent
   The major numerical tasks to be performed on each supernode are:
   \begin{description}
      \item[Factor] $A_{11} = L_{11}D_{11}L_{11}^T$;
      \item[Solve] $L_{21} = A_{21}(D_{11}L_{11}^T)^{-1}$;
      \item[Form] $S = L_{21}D_{11}L_{21}^T$ (Schur complement); and
      \item[Scatter] $S$ across other supernodes (using either multifrontal or
         supernodal techniques).
   \end{description}
   Here $L_{11}$ is unit lower triangular, and $D_{11}$ is block diagonal with
   $1\times1$ and $2\times2$ blocks. In practice, permutations are used for
   pivoting, however we omit these above for clarity of notation.
   For numerical stability, when selecting pivots the factor
   task needs to take account of the values of the entries in $A_{21}$ as well as those
   in $A_{11}$.
   For this reason, the factor and solve tasks are often combined into a
   single kernel.

   In addition to the scatter task, the key difference from an otherwise
   equivalent dense factorization is that pivots are only selected from
   within $A_{11}$. If a candidate pivot is found to be unsuitable, it is moved to a later
   supernode for elimination, with a guarantee that all pivots will be
   eliminated in the final supernode. Such pivots are said to be {\it delayed}.
   They generate additional floating-point operations and storage requirements.
   If pivots were instead chosen from $A_{21}$, much larger amounts of
   additional storage and computation would be required.
   
   With the advent of manycore processors and the growing gap between the speed
   of communication and computation, many algorithms must be rewritten to reflect
   the changing balance in resource. As the pivoting decisions must be taken
   in a serial fashion, they are highly sensitive to the latency and speed of
   any communication or bandwidth costs incurred. With current algorithms
   that take account of the entries in $A_{11}$ and $A_{21}$, all
   threads working on a supernode must endure stop-start parallelism for every
   column of the supernode. Even when running serially, performance issues are
   encountered if the entire supernode does not fit in the smallest level of
   cache.

   This paper seeks to address these issues by developing effective pivoting
   strategies that significantly reduce the amount of communication required.
   A provably stable algorithm and a heuristic algorithm are presented; we will
   refer to these algorithms as {\it compressed threshold pivoting} algorithms. The heuristic
   algorithm is faster than the provably stable alternative and it more
   accurately approximates the behaviour of traditional threshold partial
   pivoting in terms of modifications to the pivot sequence. While it
   can demonstrably fail to control the growth factor for
   some pathological constructed examples, in practice it achieves numerical
   robustness even on the most difficult practical problems.

   The rest of this paper proceeds as follows. In Section~\ref{Sec:TPP}, the
   standard threshold partial pivoting technique used in many current sparse
   symmetric indefinite codes is reviewed. Section~\ref{Sec:dense techniques} explores the
   applicability of recent work on communication-reducing pivoting for dense
   factorizations to the sparse case and reviews techniques
   that are currently used for sparse problems. The new compressed threshold pivoting
   algorithms are introduced in Section~\ref{Sec:STP}; stability and
   communication are analysed in Sections~\ref{Sec:stability}~and~\ref{Sec:comm},
   respectively. Numerical experiments are presented in Section~\ref{Sec:expts}
   and conclusions are summarized in Section~\ref{Sec:conclusions}.

   \section{Threshold partial pivoting (TPP) within a sparse direct solver}
   \label{Sec:TPP}
   \setcounter{equation}{0}
\setcounter{table}{0}
\setcounter{figure}{0}
  % \lstset{basicstyle=\ttfamily,language=pseudocode,mathescape=true, label=threshold alg}
   \lstset{language=pseudocode,mathescape=true}

   \begin{lstlisting}[float,caption=Threshold partial pivoting (TPP) algorithm,label=threshold alg]
      Input: $a(1\!:\!n,\,1\!:\!p)$ with $n \ge p$; parameters $u$ and $small$
      $nelim = 0$ // Number of eliminated variables
      $m = 2$ // Index of current pivot candidate column
      do while ( elimination still possible )
         if ( $\max( |a(m\!:\!n,\, m)| ) < small$ ) then
            permute $m$ to position $nelim+1$
            record a zero pivot; $nelim=nelim+1$// Special $1\times1$ case
            $m=m+1$; cycle // Move to next column
         end if
         Find column index $t$ of largest entry in $|a(m,\, nelim\!+\!1\!:\!m\!-\!1)|$

         // Try (t,m) as a $2\times2$ pivot
         $maxt = \max \{ |a(i,t)| :i \ge nelim+1, i \neq t, m\}$
         $maxm = \max \{ |a(i,m)| :i \ge nelim+1, i \neq t, m\}$
         if ( test_2x2 ($a$, $m$, $t$, $maxm$, $maxt$) ) then
            permute $t$ and $m$ to positions $nelim+1$ and $nelim+2$
            perform $2\times 2$ pivot
            $nelim=nelim+2$
            update $a(nelim\!+\!1\!:\!n,\, nelim\!+\!1\!:\!p)$
            $m=m+2$; cycle // Move to next column
         end if

         // Failed as $2\times2$ pivot, try as $1\times1$
         $ maxm = \max\{ |a(i,m)| : i \ge nelim+1, i \neq m$}
         if ( $|a(m,m)| \ge u*maxm$ ) then
            permute $m$ to position $nelim+1$
            perform $1\times 1$ pivot
            $nelim=nelim+1$
            update $a(nelim\!+\!1\!:\!n,\, nelim\!+\!1\!:\!p)$
            $m=m+1$; cycle // Move to next column
         end if
      end do

      // Return true if (t,m) is a good $2\times2$ pivot, false otherwise
      function test_2x2 ($a$, $t$, $m$, $maxm$, $maxt$)
         if ( $\max( |a(t,t)|, |a(t,m)|, |a(t,t)| ) < small$ ) then return false

         // Next test ensures $2\times2$ candidate is not singular and cancellation
         //  does not adversely affect the calculation of its inverse
         $detscale = 1 / \max( |a(t,t)|, |a(t,m)|, |a(t,t)| )$
         $detpiv1 = (a(t,m)*detscale) * a(t,m)$
         $detpiv0 = a(m,m)*detscale*a(t,t)$
         $detpiv = detpiv0 - detpiv1$
         if ( $| detpiv | > \max( small, |detpiv0|/2, |detpiv1|/2 )$ ) then return false

         if ( $\max(maxm, maxt) < small$ ) return true
         if ( $detpiv^{-1}*\left(detscale*
            \left | \pmat{a(m,m) & a(t,m) \\ a(t,m) & a(t,t) } \right|\right) 
             \left(\begin{array}{c}
               maxm \\
               maxt
            \end{array}\right) \le u^{-1}$ ) return true
      end function test_2x2
   \end{lstlisting}
   In this section, we recall how threshold partially pivoting (TPP) may be
   incorporated within the combined factor and solve tasks.
   
   The algorithm (which we refer to as the {\it threshold pivoting
   algorithm}) is applied to the supernodal matrix of
   (\ref{supernode}) and tries to select up to $p$ pivots from the
   first $p$ rows ($A_{11}$). The entries in the remaining rows are used when
   testing for stability.
   As pivots can only be selected from $A_{11}$, traditional partial pivoting
   is not applicable. Instead, a threshold test is employed to limit the growth
   of entries in the factors. For a $1\times1$ pivot, the test 
   on the suitability of column $q$ is 
   \begin{equation}
   |a(q,q)| \ge  u\; \max_{i>q} |a(i,q)|,
   \end{equation}
   where $a(i,j)$ are the entries of the supernodal matrix and we are assuming
   that columns $1,2,...q-1$ have already been pivoted on (see, for example, \cite{duer:86}).
   Similarly, for a $2\times 2$ pivot, the test 
   on the suitability of columns $q$ and $q+1$ is
   \begin{equation} \label{eq:2x2} 
      \left| \pmat{a(q,q) & a(q,q+1) \\ a(q,q+1) & a(q+1,q+1) }^{-1}
      \right| 
      \pmat{\max_{i>q+1} | a(i,q)| \\ \max_{i>q+1} | a(i,q+1)|}  
      \le
      \pmat{u^{-1} \\u^{-1}},
   \end{equation}
   where the absolute value notation for a matrix refers to the matrix of
   corresponding absolute values (see, for example, \cite{resc:2011}).
   The choice of the threshold parameter $u$ ($0<u\le0.5$) controls the balance
   between stability and sparsity in the
   factors, with a small $u$ leading to a potentially large growth factor but
   preserving sparsity. Observe that both tests require a scan of the
   candidate column(s), which must be up-to-date (all operations from previous
   pivots must have been applied).
   Stability of the factorization of symmetric indefinite systems was considered 
   by Ashcraft, Grimes and Lewis~\cite{agl:1998}, who showed that bounding the
   size of the entries of $L$, together with a backward stable scheme for solving
   $2 \times 2$ linear systems, suffices to show backward stability for the entire 
   process. Note that they found the widely used strategy of 
   Bunch and Kaufmann~\cite{buka:77} does not have this property. 

   Listing~\ref{threshold alg} outlines the kernel for performing
   the combined factor and solve tasks 
   that is used within our recent multifrontal solver \texttt{HSL\_MA97} \cite{hosc:2011a, hsl:2013}. 
   The description uses the notation $a(i\!:\!j,\, r\!:\!s)$
   to denote the submatrix consisting of rows $i$ to $j$ and
   columns $r$ to $s$.
   In addition to the threshold pivoting parameter $u$
   (default value $0.01$), the user-specified parameter $small$
   (default $10^{-20}$) is provided. All entries less than $small$ are treated
   as zero for the purposes of pivot selection. It is worth noting the care
   required to invert the $2\times2$ pivot in a stable fashion, whereby the
   pivot is scaled such that the largest entry is unity and the test
   $|detpiv|>\max(small, |detpiv0|/2, |detpiv1|/2)$ ensures both that the pivot
   is nonsingular and that cancellation does not occur.

   The TPP algorithm is similar to that used
   by, for example,  other HSL \cite{hsl:2013} sparse symmetric 
   indefinite solvers but has a preference for $2\times2$
   pivots over $1\times1$ pivots (see also \cite{resc:2011}). While this description is written in a
   right-looking fashion (i.e.\ the uneliminated part of the matrix is updated
   after each pivot selection), the actual implementation uses multiple levels
   of blocking, some of which use left-looking updates rather than right-looking
   ones for performance reasons (details are given in \cite{hosc:2011a}). 
   However, a key feature is that, at any given time,
   columns $nelim+1$ to $m$ must all be up-to-date, where $nelim$ is the
   number of pivots selected so far and $m$ is the index of the current
   candidate pivot column. Candidate pivots are only permuted to the front
   of the matrix once they have been accepted.

   \subsection{Parallel variants} \label{TPP:variants}
   In this paper, we will compare our new algorithms against the following
   three parallel variants of Listing~\ref{threshold alg} with different
   communication patterns. We assume that $P$ processors are used.
   
   \begin{description}
      \item[TPP outer update (TPP\_OU)] The supernodal matrix (\ref{supernode})
         is divided into block columns, each
         of width $nbi$ (with the width of the last block column
         adjusted as necessary). The supernodal factorization proceeds serially using left-looking updates
         within each block column. Upon completion of the factorization of a block column,
         a parallel right-looking update of the remaining block columns
         is performed.
      \item[TPP Variant A] The $n$ rows of the supernodal matrix $A$ (including $A_{11}$) are split equally between
         the processors. The processor that owns the current pivot row $m$ 
         determines its pairing $t$ and communicates the submatrix
         $\left(\begin{array}{cc} a(t,t) & a(t,m) \\
                                  a(t,m) & a(m,m)\end{array}\right)$ to
         each of the other processors. Each processor $k$ finds local maximum values
         $maxm_k$ and $maxt_k$ that are then reduced in parallel to find a
         global $maxm$ and $maxt$. The acceptance test and pivoting operations
         are performed locally, and the owner of the pivotal rows broadcasts
         the part of $A_{11}$ needed for the update to all other processors.
         The local updates are then performed.
      \item[TPP Variant B] The $n-p$ rows of $A_{21}$ are split equally between
         the processors. The $p$ rows of $A_{11}$ are replicated on every
         processor. Each processor independently finds the same pivots $m$ and
         $t$ and calculates local maximum values $maxm_k$ and $maxt_k$ that are then
         reduced in parallel to find a global $maxm$ and $maxt$. The acceptance
         test, pivoting operation and updates of all locally stored rows are
         then  performed without need for further communication.
   \end{description}

   We consider parallel schemes in which each processor controls
   $O(p^2)$ data and performs $O(p^3)$ operations. As such, supernodes with $n\gg p$ are
   of particular interest. In this case, each of the above variants  performs
   at least one communication per pivot, which can be hidden by at most $O(p^2)$
   operations, or less if blocking is used to exploit the cache architecture of
   individual processors. To improve on this, we need to consider ways of
   performing more pivot operations per communication.

   \section{Existing methods}
   \label{Sec:dense techniques}
   \setcounter{equation}{0}
\setcounter{table}{0}
\setcounter{figure}{0}

   In this section, we briefly review techniques that have been proposed
   to overcome the problem of pivoting in parallel, both in the dense case and in the sparse case.
   We consider their suitability for sparse indefinite systems.
   \subsection{Dense}

    In the dense case, pivots may be chosen from within $A_{21}$ as well as
   $A_{11}$. A number of different pivoting schemes have been proposed,
   including pairwise pivoting \cite{basd:60, sore:85} and parallel
   pivoting \cite{trsc:90}, and block variants thereof.
   Parallel pivoting is unstable. While pairwise pivoting is more
   stable, the growth factor is more than linear with respect to the matrix
   size \cite{grdx:2011}. Sparse versions of these algorithms are possible
   but have not been studied as they are likely to suffer from the same problems
   as the dense variants.

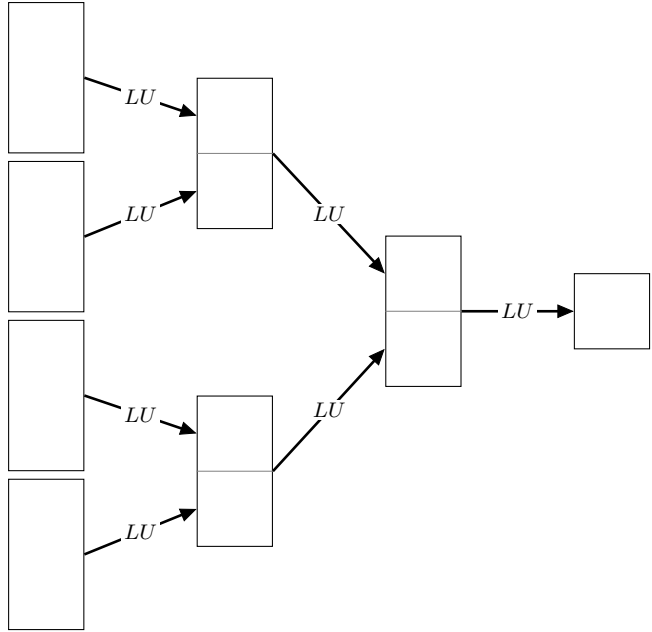
\begin{figure}
      \caption{ \label{tournament pivoting}
         Tournament pivoting on a reduction tree as used by CALU.
      }
      \vspace{3mm}
      \centering
      \begin{tikzpicture}
         \node[draw,minimum height=2cm,minimum width=1cm,inner sep=0pt] (a11) {};
         \node[draw,minimum height=2cm,minimum width=1cm,inner sep=0pt,below=0.1cm] (a21) at (a11.south) {};
         \node[draw,minimum height=2cm,minimum width=1cm,inner sep=0pt,below=0.1cm] (a31) at (a21.south) {};
         \node[draw,minimum height=2cm,minimum width=1cm,inner sep=0pt,below=0.1cm] (a41) at (a31.south) {};

         \node[draw,minimum height=2cm,minimum width=1cm,inner sep=0pt] (a12) at ($ (a11.south east) + (2,0) $) {};
         \node[draw,minimum height=2cm,minimum width=1cm,inner sep=0pt] (a22) at ($ (a31.south east) + (2,0) $) {};

         \node[draw,minimum height=2cm,minimum width=1cm,inner sep=0pt] (a13) at ($ (a12.east) + (2,-2.1) $) {};

         \node[draw,minimum height=1cm,minimum width=1cm,inner sep=0pt] (a14) at ($ (a13.east) + (2,0) $) {};

         \tikzstyle{arrow} = [
            line width=1,
            decoration={markings,mark=at position 1 with {\arrow[thin]{triangle 45}}},
            shorten >= 5.5pt,
            preaction = {decorate},
         ]

         \draw[arrow] (a11.east) -- ($ (a12.west) + (0,0.5) $) node[midway,fill=white,inner sep=1pt,font=\footnotesize] {$LU$};
         \draw[arrow] (a21.east) -- ($ (a12.west) - (0,0.5) $) node[midway,fill=white,inner sep=1pt,font=\footnotesize] {$LU$};
         \draw[arrow] (a31.east) -- ($ (a22.west) + (0,0.5) $) node[midway,fill=white,inner sep=1pt,font=\footnotesize] {$LU$};
         \draw[arrow] (a41.east) -- ($ (a22.west) - (0,0.5) $) node[midway,fill=white,inner sep=1pt,font=\footnotesize] {$LU$};

         \draw[arrow] (a12.east) -- ($ (a13.west) + (0,0.5) $) node[midway,fill=white,inner sep=1pt,font=\footnotesize] {$LU$};
         \draw[arrow] (a22.east) -- ($ (a13.west) - (0,0.5) $) node[midway,fill=white,inner sep=1pt,font=\footnotesize] {$LU$};

         \draw[arrow] (a13.east) -- (a14.west) node[midway,fill=white,inner sep=1pt,font=\footnotesize] {$LU$};

         \draw[help lines] (a12.east) -- (a12.west);
         \draw[help lines] (a22.east) -- (a22.west);
         \draw[help lines] (a13.east) -- (a13.west);

      \end{tikzpicture}
   \end{figure}

   A more radical alternative pivoting technique has recently been introduced
   by Grigori, Demmel and Xiang \cite{grdx:2011}. Their CALU algorithm uses
   {\it tournament pivoting} whereby the supernode is recursively bisected into
   sections upon which an $LU$ factorization is performed to select the best
   $p$ pivot rows. This is demonstrated in Figure~\ref{tournament pivoting}
   and motivates our approach in this paper.
   For each block, an $LU$ factorization is performed to identify the best $p$
   pivot rows within that block. These rows are then concatenated with those
   selected from a partner block and the process repeated. Once the full
   reduction tree has been evaluated, the selected pivots are used for the
   factorization of the supernode. While this technique provides
   weaker guarantees upon growth than traditional partial pivoting, it is no
   worse than partial pivoting applied at a block level. An analysis is
   presented in \cite{grdx:2011} showing that this algorithm performs an
   optimal amount of communication that is asymptotically less than Gaussian
   elimination with partial pivoting. Hence, it is faster on platforms where
   communication is expensive. Furthermore, the method has been shown to be
   stable in practice. However, as pivots are selected from within
   $A_{21}$, it is not applicable to the sparse case.

   Other approaches seek to avoid the need for pivoting altogether. Becker,
   Baboulin and Dongarra \cite{bebd:2011} use a randomizing scaling with
   structure similar
   to a fast-Fourier transform to homogenise the matrix such that the values
   are sufficiently uniformly distributed that pivoting is almost surely not
   required. This technique of using {\it recursive butterfly matrices} 
   cannot be applied to the sparse case as they lead to complete fill-in of 
   the scaled matrix $SAS$.

   \subsection{Sparse}

   Since searching candidate columns is expensive, an obvious remedy is to
   restrict the search and possibly risk sacrificing some stability.
   In a {\it restricted pivoting} approach,
   the check for large entries is restricted to examining only the  values 
   of the entries in $A_{11}$,
   completely ignoring those in $A_{21}$.
   As it may not be possible to find a
   suitable pivot from within the diagonal block, the sparse direct solver PARDISO \cite{scga:2006}
   uses restricted pivoting with a pivot
   perturbation strategy (called {\it static pivoting}), which is similar to that
   employed in the unsymmetric case in  \cite{lide:98}. Static pivoting
   allows no pivots to be delayed, thereby greatly simplifying the coding
   compared with a direct solver that does permit delayed pivots, as well as
   limiting the fill in the factors and the operations required to compute them.
   However, since the factorization may not be accurate, it is often necessary
   to perform a number of steps of  refinement to try and recovery accuracy but this is not
   guaranteed to be successful. However, as discussed by Schenk, W\"achter and
   Hagemann in \cite{scwh:2007},
   reliability may be improved by using a matching-based ordering. This aims to
   bring large entries close to the diagonal, with the hope that they will make
   suitable block pivot candidates (see also Duff and Pralet~\cite{dupr:2007}). 
   While this extends the class of problems that can be solved,
   there are some particularly tough indefinite linear systems that we have been unable to
   solve successfully using this approach (see \cite{hosc:2013}).
   Furthermore, some applications, particularly those arising from
   optimization, require an accurate measure of the matrix inertia that can be
   difficult to obtain when static pivoting is used.

   Kim and Eijkhout \cite{kiei:2012} present a {\it try-it-and-see} approach to
   pivoting in their code for $hp$-adaptive finite element problems. Restricted
   pivoting is used at each node of the assembly tree, but an a posteriori check
   is made for growth below the diagonal block. Should excessive growth be
   detected, the supernodal matrix is reconstructed from the 
   contribution blocks of its child nodes and pivots are delayed to its parent node.
 
   \section{Compressed threshold pivoting}
   \label{Sec:STP}
\setcounter{equation}{0}
\setcounter{table}{0}
\setcounter{figure}{0}
   As stated at the end of Section~\ref{Sec:TPP}, we need to reduce the number
   of communications per pivot. For $n\gg p$, we propose the
   construction of a small representative matrix that can be used to make
   pivoting decisions without the need to work with the full 
   supernodal matrix (and
   hence avoids the need to communicate with other processors).
   Henceforth we shall refer to this representative matrix as the
   {\it compressed matrix} $C$. For the purposes of this paper, we shall
   consider only the case where $C$ is $p\times p$, but there is no
   requirement that this is so.
   The application of (potentially modified) threshold partial
   pivoting to the small trapezoidal matrix comprising the diagonal block $A_{11}$
   of (\ref{supernode}) and the matrix $C$ as below is considered.
\newpage
   \begin{center}
      \begin{tikzpicture}
         \draw (0,1.1) -- (0,2.2) -- (1,1.1) -- cycle;
         \node at (0.35, 1.5) {$A_{11}$};
         \draw (0,0) -- (1,0) -- (1,1) -- (0,1) -- cycle;
         \node at (0.5,0.5) {$C$};
      \end{tikzpicture}
      \vspace*{-1.75cm}\begin{eqnarray}\label{A+C}\end{eqnarray}\vspace*{0.4cm}
   \end{center}

   \begin{sidewaysfigure}
      \caption{ \label{subset tp fig}
         Compressed threshold pivoting
      }
      \vspace{3mm}
      \begin{tikzpicture}
         \matrix {
         \node[text width=3.6cm]{\begin{center} (a)\\ Analysis of $A_{21}$\end{center}}; &
         \node[text width=3.6cm]{\begin{center} (b)\\ Construction of compressed matrix $C$\end{center}}; &
         \node[text width=4.0cm]{\begin{center} (c)\\ Factor\\ $\left(\begin{array}{c}A_{11}\\C\end{array}\right)=(P\hat{L})D(P\hat{L})^T$\end{center}}; &
         \node[text width=4.0cm]{\begin{center} (d)\\ Apply permutation\\ $P$ to $A_{21}$\end{center}};&
         \node[text width=4.0cm]{\begin{center} (e)\\ Solve \\$D_{11}L_{11}^{T}L_{21}=PA_{21}$\end{center}}; \\
            
            % (a)
            \draw (-0.5,0) -- (0.5,0) -- (-0.5,1) -- cycle; 
            \draw (0.5,0) -- (0.5,-4) -- (-0.5,-4) -- (-0.5,0);
            \node[left, inner sep=0] at (-0.5,0.94) {\scriptsize x};
            \node[left, inner sep=0] at (-0.5,0.79) {\scriptsize x};
            \node[left, inner sep=0] at (-0.5,0.64) {\scriptsize x};
            \node[left, inner sep=0] at (-0.5,0.49) {\scriptsize x};
            \node[left, inner sep=0] at (-0.5,0.34) {\scriptsize x};
            \node[left, inner sep=0] at (-0.5,0.19) {\scriptsize x};
            \node[left, inner sep=0] at (-0.5,0.04) {\scriptsize x};
            \node[left, inner sep=0] (a1) at (-0.5,-0.22) {\scriptsize x};
            \node[left, inner sep=0] (a2) at (-0.5,-0.72) {\scriptsize x};
            \node[left, inner sep=0] (a3) at (-0.5,-0.87) {\scriptsize x};
            \node[left, inner sep=0] (a4) at (-0.5,-1.42) {\scriptsize x};
            \node[left, inner sep=0] (a5) at (-0.5,-2.42) {\scriptsize x};
            \node[left, inner sep=0] (a6) at (-0.5,-3.22) {\scriptsize x};
            \node[left, inner sep=0] (a7) at (-0.5,-3.92) {\scriptsize x};
            \draw (-0.43,-0.21) circle(.07);
            \draw (-0.29,-3.91) circle(.07);
            \draw (-0.15,-0.71) circle(.07);
            \draw (-0.01,-1.41) circle(.07);
            \draw (0.13,-3.21) circle(.07);
            \draw (0.27,-0.86) circle(.07);
            \draw (0.41,-2.41) circle(.07);
            \draw[line width=1] (0.6,1.0) -- (0.7,1.0) -- (0.7,0.0) -- (0.6,0.0);
            \node[inner sep=0] (a0) at (0.7,0.5) {};
            \node at (-0.1,0.3) {$A_{11}$};
            \node at (0,-2.0) {$A_{21}$};
            \node (aA21) at (0.5,-2.0) {};
            &
            % (b)
            \draw (-0.5,0) -- (0.5,0) -- (-0.5,1) -- cycle; 
            \draw (0.5,0) -- (0.5,-1) -- (-0.5,-1) -- (-0.5,0);
            \node[left, inner sep=0] at (-0.5,0.94) {\scriptsize x};
            \node[left, inner sep=0] at (-0.5,0.79) {\scriptsize x};
            \node[left, inner sep=0] at (-0.5,0.64) {\scriptsize x};
            \node[left, inner sep=0] at (-0.5,0.49) {\scriptsize x};
            \node[left, inner sep=0] at (-0.5,0.34) {\scriptsize x};
            \node[left, inner sep=0] at (-0.5,0.19) {\scriptsize x};
            \node[left, inner sep=0] at (-0.5,0.04) {\scriptsize x};
            \node[left, inner sep=0] (b1) at (-0.5,-0.09) {\scriptsize x};
            \node[left, inner sep=0] (b2) at (-0.5,-0.23) {\scriptsize x};
            \node[left, inner sep=0] (b3) at (-0.5,-0.37) {\scriptsize x};
            \node[left, inner sep=0] (b4) at (-0.5,-0.51) {\scriptsize x};
            \node[left, inner sep=0] (b5) at (-0.5,-0.65) {\scriptsize x};
            \node[left, inner sep=0] (b6) at (-0.5,-0.79) {\scriptsize x};
            \node[left, inner sep=0] (b7) at (-0.5,-0.91) {\scriptsize x};
            \node[inner sep=0] (b0) at (-0.5,0.5) {};
            \node at (-0.1,0.3) {$A_{11}$};
            \node at (0,-0.5) {$C$};
            \node[inner sep=0] (bmid) at (0.6,0.0) {};
            &
            % (c)
            \draw (-0.5,0) -- (0.5,0) -- (-0.5,1) -- cycle; 
            \node at (-0.1,0.3) {$L_{11}$};
            \draw (0.5,0) -- (0.5,-1) -- (-0.5,-1) -- (-0.5,0);
            \node at (0,-0.5) {$\hat{C}$};
            \draw (-1.6,0) -- (-0.6,0) -- (-0.6,1) -- (-1.6,1) -- cycle;
            \node at (-1.1,0.5) {$P$};
            \draw (1.6,0) -- (0.6,0) -- (0.6,1) -- (1.6,1) -- cycle;
            \node at (1.1,0.5) {$D_{11}$};
            \draw[line width=1] (1.7,1.0) -- (1.8,1.0) -- (1.8,0.0) -- (1.7,0.0);
            \draw[line width=1] (-1.7,1.0) -- (-1.8,1.0) -- (-1.8,0.0) -- (-1.7,0.0);
            \node[inner sep=0] (cmid) at (-1.8,0.0) {};
            \node[inner sep=0] (cPbottom) at (-1.1,0.0) {};
            \node[inner sep=0] (c0) at (1.8,0.5) {};
            &
            % (d)
            \draw (-0.5,0) -- (0.5,0) -- (-0.5,1) -- cycle; 
            \node at (-0.1,0.3) {$L_{11}$};
            \draw (0.5,0) -- (0.5,-4) -- (-0.5,-4) -- (-0.5,0);
            \node at (0,-2.0) {$PA_{21}$};
            \draw (-1.6,0) -- (-0.6,0) -- (-0.6,1) -- (-1.6,1) -- cycle;
            \node at (-1.1,0.5) {$P$};
            \draw (1.6,0) -- (0.6,0) -- (0.6,1) -- (1.6,1) -- cycle;
            \node at (1.1,0.5) {$D_{11}$};
            \draw[line width=1] (1.7,1.0) -- (1.8,1.0) -- (1.8,0.0) -- (1.7,0.0);
            \draw[line width=1] (-1.7,1.0) -- (-1.8,1.0) -- (-1.8,0.0) -- (-1.7,0.0);
            \node[inner sep=0] (d0left) at (-1.8,0.5) {};
            \node[inner sep=0] (d0right) at (1.8,0.5) {};
            \node (dmid) at (-0.5,-2.0) {};
            \node (dA21) at (0.5,-2.0) {};
            &
            % (e)
            \draw (-0.5,0) -- (0.5,0) -- (-0.5,1) -- cycle; 
            \node at (-0.1,0.3) {$L_{11}$};
            \draw (0.5,0) -- (0.5,-4) -- (-0.5,-4) -- (-0.5,0);
            \node at (0,-2.0) {$L_{21}$};
            \draw (-1.6,0) -- (-0.6,0) -- (-0.6,1) -- (-1.6,1) -- cycle;
            \node at (-1.1,0.5) {$P$};
            \draw (1.6,0) -- (0.6,0) -- (0.6,1) -- (1.6,1) -- cycle;
            \node at (1.1,0.5) {$D_{11}$};
            \node[inner sep=0] (emid) at (-0.5,-2.0) {};
            \\
         };
         \draw ($ (a1) + (1.1,0) $) edge [->,bend left=15] (b1);
         \draw ($ (a2) + (1.1,0) $) edge [->,bend left=20] (b2);
         \draw ($ (a3) + (1.1,0) $) edge [->,bend left=20] (b3);
         \draw ($ (a4) + (1.1,0) $) edge [->,bend left=20] (b4);
         \draw ($ (a5) + (1.1,0) $) edge [->,bend left] (b5);
         \draw ($ (a6) + (1.1,0) $) edge [->,bend left] (b6);
         \draw ($ (a7) + (1.1,0) $) edge [->,bend left] (b7);
         \draw (a0) edge [->,line width=1] ($(b0)+(-0.2,0.0)$);
         \draw[vecArrow] (bmid) -- (cmid) node[midway,below] {Factor};
         \draw (c0) edge [->,line width=1] (d0left);
         \draw[line width=1] (aA21) to[in=180,out=320] ($ (dmid) + (-1.5,-0.05) $);
         \draw[line width=1] (cPbottom) to[in=180,out=270] ($ (dmid) + (-1.5,0.05) $);
         \draw[vecArrow] ($ (dmid) + (-1.5,0) $) -- (dmid) node[midway,below] {Permute};
         \draw[line width=1] (dA21) to[in=180,out=320] ($ (emid) + (-1.0,-0.05) $);
         \draw[line width=1] (d0right) to[in=180,out=320] ($ (emid) + (-1.0,0.05) $);
         \draw[vecArrow] ($ (emid) + (-1.0,0) $) -- (emid) node[midway,below] {Solve};
      \end{tikzpicture}
   \end{sidewaysfigure}
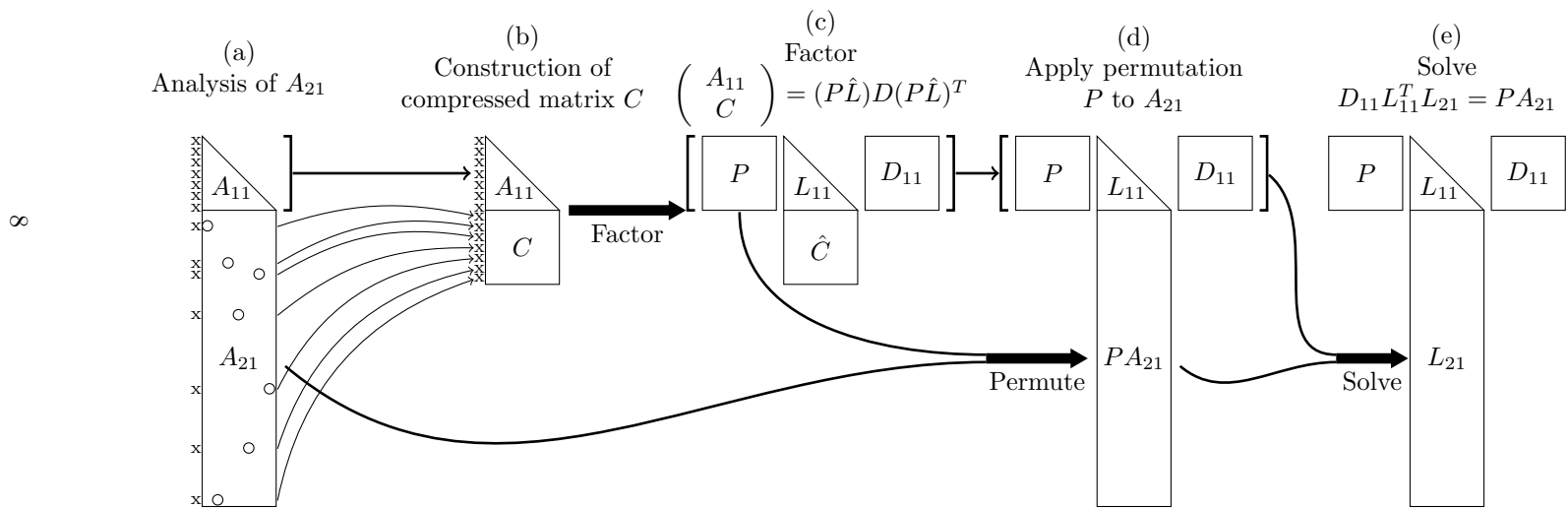
   
   \noindent
   The factorization of (\ref{A+C}) establishes a permutation and factors $L_{11}$
   and $D_{11}$ that
   can be applied to $A_{21}$ (in parallel) without the need to perform any further
   pivoting. This process is summarized in Figure~\ref{subset tp fig}.

   We present two methods of constructing (and updating) the compressed matrix.
   The {\it strict} method is numerically
   stable, however it can be too pessimistic (that is, lead to a large number
   of delayed pivots that TPP would have selected) so we also introduce the 
   {\it relaxed} method that can be unstable for some matrices but, in practice,
   if used with an appropriate scaling and ordering, will be shown to be stable (see
   Section~\ref{real-world results}).

   Both methods construct their compressed matrix in parallel using a tree
   reduction similar to that previously described for the CALU algorithm. The $LU$
   reduction operation of Figure~\ref{tournament pivoting}
   is replaced by the techniques that we describe in the
   next two sections.

   \subsection{Strict compressed pivoting}

   In the strict method, we first partition the rows of $A_{21}$ into sets
   corresponding to the column in which their entry of largest absolute value lies.
   Thus all the indices of the rows in $A_{21}$ that have their
   maximum entry in column $j$ belong to  the set  $J_j$ given by 
   $$
      J_j=\{i : j=\arg\max_k |a(i, k)|, i>p\}.
   $$
   Ties are resolved in favour of the lowest value of $k$.
   The compressed matrix $C=\{c(j,\,k)\}$ is constructed row-by-row from these
   sets:
   \begin{equation}\label{eq:scp}
      c(j,k) = \left\{\begin{array}{ll}
            \max_{i\in J_j} |a(i,k)|, & J_j\ne\phi, \\
            0, & \mathrm{otherwise}.\\
         \end{array}\right.
   \end{equation}
   That is, row $j$ of $C$ is determined by the columnwise maxima of $J_j$. This
   is illustrated in Figure~\ref{strict eg}.
   
   \begin{figure}[htbp]
      \caption{ \label{strict eg}
         Example of strict compressed matrix construction.
         Observe that $J_1=\phi$. Bold is used to indicate row maxima of
         $A_{21}$ and column maxima of $A_{J_2}$ and $A_{J_3}$.
      }
      \vspace{3mm}
      \centering
      \begin{tikzpicture}
         \node (a21) {$A_{21}=\left(\begin{array}{rrr}
                  1 & \bf 10 & 10 \\
                  2 & 3      & \bf 4 \\
                    & \bf 10 & -3 \\
                  4 & \bf -5 & 4 \\
                    & -6     & \bf 8 \\
               \end{array}\right)$};
         \node[right] (aj2) at ($ (a21.east) +(0.5,0.50) $) {
            $A_{J_2}=\left(\begin{array}{rrr}
            1 & \bf 10 & \bf 10 \\
            & 10 & -3 \\
            \bf 4 & -5 & 4
         \end{array}\right)$};
         \node[below] (aj3) at (aj2.south) {$A_{J_3}=\left(\begin{array}{rrr}
            \bf 2 & 3 & 4 \\
            & \bf -6 & \bf 8 \\
            \end{array}\right)$
            };
         \node[right] at ($ (aj2.south east) +(0.5,0.2)$) {$C=\left(\begin{array}{rrr}
            0 & 0 & 0 \\
            4 & 10 & 10 \\
            2 & 6 & 8 \\
            \end{array}\right)$
            };
      \end{tikzpicture}
   \end{figure}

   The factorization proceeds as per threshold partial pivoting applied to
   (\ref{A+C}). However, the steps involving the application of pivots and
   the updating of the trailing submatrix are modified for $C$ (but not for
   $A_{11}$). For threshold partial pivoting we would use the following update
   formulae for a $1\times1$ pivot:
   \begin{eqnarray*}
      \hat{c}(:,\,k) &=& c(:,\,k)/a(k,\,k) \\
      \hat{c}(:,\,k\!+\!1\!:\!p) &=& c(:,\,k\!+\!1\!:\!p) - c(:,\,k) a(k,\,k) a(k\!+\!1\!:\!p,\,k)^T.
   \end{eqnarray*}
   \vspace*{-1.7cm}\begin{eqnarray}\label{eq:pivot1x1}\end{eqnarray}\vspace*{0.0cm}\\
   We want row $j$ of $C$ to represent the worst possible growth in the
   rows $J_j$. Using similar ideas to those involved in the $2\times2$ pivot
   test (\ref{eq:2x2}), we modify (\ref{eq:pivot1x1}) to use absolute
   values throughout (exploiting the fact that all elements of $C$ are positive):
   \begin{eqnarray*}
      \hat{c}(:,\,k) &=& c(:,\,k)/|a(k,\,k)| \\
      \hat{c}(:,\,k\!+\!1\!:\!p) &=& c(:,\,k\!+\!1\!:\!p) + c(:,\,k)\, |a(k,\,k)|\, |a(k\!+\!1\!:\!p,\,k)|^T.
   \end{eqnarray*}
   \vspace*{-1.7cm}\begin{eqnarray}\label{eq:worst case upd}\end{eqnarray}\vspace*{0.0cm}\\
   A proof of backwards stability is given in Section~\ref{Sec:stability}.
   
   We observe that we experimented with using a single row to represent the entire
   $A_{21}$ matrix (based on the set $J=\{1,2,\ldots,n-p\}$, that is, $C$ is
   constructed by taking the entry of largest absolute value in each row
   of  $A_{21}$), but found this led to the rejection
   of almost all pivots beyond the first few. By using multiple rows, the
   over estimation of the growth can be controlled.

   \subsection{Relaxed compressed pivoting}

   For the relaxed method, we make
   the experimental observation that the column maxima often remain in the
   same position as the factorization progresses. Even in those cases where the
   locations of the maxima vary, the value at the old location is often close to
   that at the new.
   
   For each column of the matrix (\ref{supernode}), we include in the
   compressed matrix a row that contains the entry of largest absolute value in
   that column
   at the start of the supernode factorization. The hope is that this row carries
   sufficient information to reject unstable pivots using the standard tests and
   normal update formulae (as given in (\ref{eq:pivot1x1})).
   However, the method risks using
   unstable pivots rather than rejecting acceptable ones.

   The algorithm for constructing $C$ is as follows. First mark all rows of
   $A_{21}$ as unflagged. Then, for each column $j$ ($1 \le j \le p$), find the
   the entry of largest
   magnitude in $A_{21}$ outwith an already flagged row. Flag the corresponding
   row, include it in $C$, and continue to the next column. The results of this
   algorithm are demonstrated in Figure~\ref{relaxed eg}.

   \begin{figure}[htbp]
      \caption{ \label{relaxed eg}
         Example of relaxed compressed matrix construction.
         Bold is used to indicate the entries of $A_{21}$ used in the selection
         of rows for inclusion in $C$.
      }
      \vspace{3mm}
      \centering
      \begin{tikzpicture}
         \node (a21) {$A_{21}=\left(\begin{array}{rrr}
                  1 & \bf 10 & 10 \\
                  2 & 3      & 4 \\
                    & 10     & -3 \\
              \bf 4 & -5     & 4 \\
                    & -6     & \bf 8 \\
               \end{array}\right)$};
         \node[right] at ($ (a21.east) +(1.0,0.0)$) {$C=\left(\begin{array}{rrr}
            4 & -5 &  4 \\
            1 & 10 & 10 \\
              & -6 &  8 \\
            \end{array}\right)$
            };
      \end{tikzpicture}
   \end{figure}

   Note that by insisting on flagging a different row for each column, $p$ rows
   are always marked, and ordering affects tie-breaking (for example, in
   Figure~\ref{relaxed eg} encountering row 3 before row 1 would have
   resulted in a different $C$). The choice to add the
   complication of flagging is based on practical experience.
   Otherwise, if for several columns the  entry of largest absolute value
   occurred in the
   same row, the total number of included rows could be significantly fewer
   than $p$.
   Experiments showed that in this case the resulting factorization
   was less stable for a number of problems tested.
   
   %Our numerical results will show that, in many practical examples, this
   %method achieve convergence on more problems than using restricted pivoting
   %(equivalent to taking $C=0$).

   \section{Stability analysis} \label{Sec:stability}
   Following the stability analysis presented by Ashcraft, Grimes
   and Lewis \cite{agl:1998}, we define the partially factorized supernodal matrix $A^{(q)}$
   as that formed after $q$ eliminations and their updates have been applied to
   (\ref{supernode}), with $A^{(0)}=A$. Let
   $$\mu_q=\max_{i,j} |a^{(q)}(i, j)|$$
   be the maximum absolute value of an entry of $A^{(q)}$. We proceed to analyse
   the application of the algorithm given in Listing~\ref{threshold alg} to
   (\ref{A+C}). The analysis
   assumes that the pivotal columns $t$ and $m$ have been permuted to positions
   $q$ and $q+1$, that is to be the first columns in the uneliminated part of
   the matrix.
   We seek to demonstrate that the entries of the factor $L$
   are bounded, and that growth in $A^{(q)}$ (and hence $D$) is limited:
   \begin{equation}\label{eq:bnd}
       \mu_{q+1}<\mu_{q}\left(1+u^{-1}\right).
   \end{equation}

   We first note that (\ref{eq:bnd}) does not hold for relaxed compressed pivoting.
   To illustrate this, consider the following factorization, where $\epsilon$ is small,
   $$
      A=\left(\begin{array}{cc}
         1 & -1 \\
         -1 & 2 \\
         u^{-1} \\
         & u^{-1} \\
         \hline
         u^{-1}-\epsilon & u^{-1}-\epsilon \\
      \end{array}\right)
      \qquad\Rightarrow\qquad
      L=\left(\begin{array}{cc}
         1 \\
         -1 & 1 \\
         u^{-1} & u^{-1}\\
                & u^{-1} \\
         \hline
         u^{-1}-\epsilon & 2(u^{-1}-\epsilon) \\
      \end{array}\right).
   $$
   The entries below the line are not included in the compressed matrix, so are
   not tested for stability. Observe that after $p$ steps, we can have an
   entry of $L$ that is close to $pu^{-1}$. This means that $L$ can
   be effectively unbounded.

   We next analyse strict compressed pivoting. We have the following
   bound on the entries of $A_{21}^{(q)}$.
  \begin{lemma} \label{SCP1}
  Let the compressed matrix $C$ be defined by (\ref{eq:scp}) and
   let $C^{(q)}$ to be the matrix $C$
   after $q$ eliminations ($C^{(0)}=C$). Then for $q \ge 0$,
   $$|a^{(q)}(i,k)| \le c^{(q)}(j,k)\;\;\; \mbox{   for all   } \;\;\;i\in J_j,1\le j,k\le p.$$
\end{lemma}

\noindent {\bf Proof:} 
   This
   is by construction for $q=0$, and induction on the update equations
   (\ref{eq:worst case upd}) for $q>0$ (for notational convenience here
   and elsewhere we drop
   the superscript on the reduced matrix at step $q$, but not for other steps):
   \begin{eqnarray*}
      a^{(q+1)}(i,k) & = & a(i,k) - \frac{a(i,q)a(k,q)}{a(q,q)} \\
      \Rightarrow\left|a^{(q+1)}(i,k)\right| & \le &
         c(j,k) + \frac{c(j,q)|a(k,q)|}{|a(q,q)|} = c^{(q+1)}(j,k)
   \end{eqnarray*}
   \hfill $\Box$
  %% Similarly, $|l(i,q)|\le c(i,q)$. 

      \vspace*{0.5cm}
   This lemma is used to prove strict compressed pivoting is backwards stable.
   \begin{theorem}\label{SCP2}
    For strict compresssed pivoting the bound (\ref{eq:bnd}) holds and the
    entries of $L$ are bounded above by $u^{-1}$.
   \end{theorem}

   \noindent {\bf Proof:} 
   We  proceed as in \cite{agl:1998} and consider $1 \times 1$ and $2 \times 2$ pivots separately.

   \vspace*{0.5cm}

   \noindent
   \textbf{$\mathbf{1\times1}$ pivots.} Define
   $$\gamma_q=\max\left(\max_{q<i\le p}|a(i,q)|, \max_{1\le i\le p}c(i,q)\right).$$
   Then if $a(q,q)$ is a $1\times1$ pivot it satisfies $a(q,q) \ge u\gamma_q$.
   The entries of $A^{(q)}$ are given by
   \begin{eqnarray*}
      a^{(q+1)}(i,k) & = & a(i,k) - \frac{a(i,q)a(k,q)}{a(q,q)}. \\
   \end{eqnarray*}
   Using Lemma~\ref{SCP1}, these may be bounded by
   \begin{eqnarray*}
      |a^{(q+1)}(i,k)| 
      & \le & \left\{\begin{array}{l@{\qquad}l}
         
      \displaystyle |a(i,k)| + \frac{|a(i,q)|\;|a(k,q)|}{|a(q,q)|} & i\le p\\[1em]
         \displaystyle |a(i,k)| + \frac{|c(i,q)|\;|a(k,q)|}{|a(q,q)|} & i>p\\
         \end{array}\right. \\
      & \le & |a(i,k)| + \frac{\gamma_q^2}{|a(q,q)|}. \\
   \end{eqnarray*}
   Taking maximums over $i, k$ yields
   \begin{eqnarray*}
      \mu_{q+1} & \le & \mu_q + \frac{\gamma_q^2}{|a(q,q)|}\\
      & \le & \mu_q + \gamma_qu^{-1}   \le    \mu_q\left(1+u^{-1}\right),
   \end{eqnarray*}
   and $A^{(q+1)}$ has its growth bounded. Similarly, the entries of $L$
   are bounded,
   \begin{eqnarray*}
      |l(i,q)| &=& \frac{|a(i,q)|}{|a(q,q)|} \le u^{-1}\frac{a(i,q)}{\gamma_q} \le u^{-1}.
   \end{eqnarray*}

   \vspace*{0.5cm}

   \noindent
   \textbf{$\mathbf{2\times2}$ pivots.} Define
   \begin{eqnarray*}
      \gamma_q&=&\max\left(\max_{q+1<i\le p}a(i,q), \max_{1\le i\le p}c(i,q)\right),\\
      \gamma_{q+1}&=&\max\left(\max_{q+1<i\le p}a(i,q+1), \max_{1\le i\le p}c(i,q+1)\right),\\
      D_q & = & \left(\begin{array}{cc}
            a(q,q)   & a(q,q+1) \\
            a(q,q+1) & a(q+1,q+1)
      \end{array}\right).
   \end{eqnarray*}
   Then the nonsingular pivot $D_q$ satisfies
   \begin{eqnarray*}
      \left|D_q^{-1}\right|\left(\begin{array}{c}
            \gamma_q \\ \gamma_{q+1}
      \end{array}\right)
      \le
      \left(\begin{array}{c}
            u^{-1} \\
            u^{-1}
      \end{array}\right).
   \end{eqnarray*}
   The entries of $A^{(q+2)}$ are given by
   \begin{eqnarray*}
      a^{(q+2)}(i,k) & = & a(i,k) + \left(\begin{array}{cc}
            a(k,q) &
            a(k,q+1)
      \end{array}\right)D_q^{-1}\left(\begin{array}{c}
            a(i,q) \\
            a(i,q+1)
      \end{array}\right).
   \end{eqnarray*}
   Again, using Lemma~\ref{SCP1}, these may be bounded by
   \begin{eqnarray*}
      |a^{(q+2)}(i,k)| & \le & \left\{\begin{array}{l@{\qquad}l}
         |a(i,k)| + \left(\begin{array}{cc}
            |a(k,q)| &
            |a(k,q+1)|
         \end{array}\right)\left|D_q^{-1}\right|\left(\begin{array}{c}
            a(i,q) \\
            a(i,q+1)
      \end{array}\right) & i\le p\\[1em]
         |a(i,k)| + \left(\begin{array}{cc}
            |a(k,q)| &
            |a(k,q+1)|
         \end{array}\right)\left|D_q^{-1}\right|\left(\begin{array}{c}
            c(i,q) \\
            c(i,q+1)
      \end{array}\right) & i>p
      \end{array}\right.\\[0.5em] 
      & \le & |a(i,k)| + \left(\begin{array}{cc}
            |a(k,q)| &
            |a(k,q+1)|
      \end{array}\right)\left|D_q^{-1}\right|\left(\begin{array}{c}
            \gamma_q \\
            \gamma_{q+1}
      \end{array}\right) \\
      & \le & |a(i,k)| + \left(\begin{array}{cc}
            |a(k,q)| &
            |a(k,q+1)|
      \end{array}\right)\left(\begin{array}{c}
            u^{-1} \\
            u^{-1}
      \end{array}\right). \\
   \end{eqnarray*}
   Taking maximums over $i,k$ yields
   \begin{eqnarray*}
      \mu_{q+2} & \le &
         \mu_q + u^{-1} \left(|a(k,q)| + |a(k,q+1)|\right) \\
      & \le & \mu_q(1+2u^{-1}).
   \end{eqnarray*}
   Likewise the entries of $L$ are bounded,
   \begin{eqnarray*}
      \left(\begin{array}{cc}
            l(i,q) & l(i,q+1)
         \end{array}\right) & = & \left(\begin{array}{cc}
            a(i,q) &
            a(i,q+1)
         \end{array}\right)D_q^{-T} \\
      \left(\begin{array}{cc}
            |l(i,q)| & |l(i,q+1)|
         \end{array}\right) & \le & \left(\begin{array}{cc}
            \gamma_q &
            \gamma_{q+1}
         \end{array}\right)\left|D_q^{-T}\right| \\
      & \le & \left(\begin{array}{cc}
            u^{-1} & u^{-1}
         \end{array}\right).
   \end{eqnarray*}
   \hfill $\Box$

   \section{Communication analysis} \label{Sec:comm}
   \setcounter{equation}{0}
\setcounter{table}{0}
\setcounter{figure}{0}
   Appendix~\ref{appendix:comm anal} presents an analysis of communication
   costs using a model of a parallel machine. We count the total number of
   operations and total amount of bandwidth required. We also derive a count
   of the number of messages required as a count the number of communication
   latencies that are necessarily incurred if all other operations take zero
   time. We assume that each pivot is accepted as soon as it is encountered (and
   hence permutations are not applied).

   Table~\ref{comm:actual} gives exact results for
   factorizing an $n\times p$ supernodal matrix on $P$ processors. Under the assumption
   that  $P=O(n)$, 
   Table~\ref{comm:order} summarises the results using order notation. The 
   operation counts are given in terms of the number of additional operations
   above those required for traditional threshold partial pivoting,
   $$
      \mathrm{TPP}_{\mathrm{ops}}(n,p) = \frac{29}{6}p -\frac{3}{4}p^2 -\frac{1}{3}p^3 + 2np + \frac{1}{2}np^2.
   $$
   We observe that restricted pivoting provides a theoretical bound on the
   number and amount of communication for a parallel code of the type considered
   here, as it performs no communication for pivoting. The main gain from using
   compressed pivoting compared with Variants A and B is the factor $p$ reduction in the number of messages
   sent, at the cost of using five times more bandwidth. In terms of operations,
   only Variant B significantly increases the (asymptotic) operation count
   above that of threshold partial pivoting, although strict compressed pivoting does
   add an extra $Pp^2$ term. However, these increased operation counts are
   somewhat misleading as they are spread across $P$ processors.
   
   We remark that, for a sparse direct solver, the costs must be summed over all the
   supernodes. If there are delayed pivots at a supernode then $n$ and $p$ will increase
   beyond that predicted by the analyse phase of the solver, leading
   to the need to perform permutations,
   the repeated testing of candidate pivots, additional operations and communication 
   as well as denser factors. For strict compressed pivoting,
   the number of delayed pivots is usually greater than for relaxed compressed pivoting
   and can be prohibitive; this
   is illustrated by the results in Section~\ref{real-world results}.

   \begin{table}[htbp]
      \caption{ \label{comm:actual}
         Parallel communication analysis results
      }
      \vspace{3mm}
      \centering
      \small
      \begin{tabular}{l|l|c|c}
         & Operations & Messages & Bandwidth \\
         \hline
         TPP Variant A &
            $\mathrm{TPP}_{\mathrm{ops}}(n,p)$ &
            $p+\frac{1}{2}p\log_2P$ &
         $-\frac{1}{2}p + \frac{1}{2}Pp(p+2)$ \T \\ 
         TPP Variant B &
            $\mathrm{TPP}_{\mathrm{ops}}(n,p) + (P-1)(\frac{29}{6}p + \frac{5}{4}p^2 +\frac{1}{6}p^3)$ &
            $1+\frac{1}{2}p\log_2P$ &
            $-\frac{1}{2}p(p+3)+\frac{1}{2}Pp(p+5)$ \T \\
         Strict Compressed & 
            $\mathrm{TPP}_{\mathrm{ops}}(n,p) + \frac{1}{2}p((p-1)p+3) + n(2p-1) + Pp^2$ &
            $1 + \log_2 P$ &
            $-\frac{1}{2}p(5p+1) + \frac{1}{2}Pp(5p+1)$ \T \\ 
         Relaxed Compressed & 
            $\mathrm{TPP}_{\mathrm{ops}}(n,p) + \frac{1}{2}p((p+2)p-2) + (n+P)p$ &
            $1 + \log_2 P$ &
            $-\frac{1}{2}p(5p+1) + \frac{1}{2}Pp(5p+1)$ \T \\
         Restricted & 
            $\mathrm{TPP}_{\mathrm{ops}}(n,p) - p(n-p)$ &
            $1$ &
            $-\frac{1}{2}p(p+1) + \frac{1}{2}Pp(p+1)$ \T \\
         \hline
      \end{tabular}
   \end{table}

   \begin{table}[htbp]
      \caption{ \label{comm:order}
         Summary of parallel communication analysis results (assuming  $P=O(n)$).
      }
      \vspace{3mm}
      \centering
      \small
      \begin{tabular}{l|c|c|c}
         & Operations & Messages & Bandwidth \\
         \hline
         TPP Variant A &
            $O(np^2)$ &
            $O(p\log n)$ &
            $O(np^2)$ \T \\
         TPP Variant B &
            $O(np^3)$ &
            $O(p\log n)$ &
            $O(np^2)$ \T \\
         Relaxed Compressed & 
            $O(np^2)$ &
            $O(\log n)$ &
            $O(np^2)$ \T \\
         Strict Compressed & 
            $O(np^2)$ &
            $O(\log n)$ &
            $O(np^2)$ \T \\
         Restricted & 
            $O(np^2)$ &
            $O(1)$ &
            $O(np^2)$ \T \\
         \hline
      \end{tabular}
   \end{table}

   \section{Numerical experiments} \label{Sec:expts}
\setcounter{equation}{0}
\setcounter{table}{0}
\setcounter{figure}{0}
   Because of the complexity involved in efficiently implementing all our algorithms and
   measuring timing performance within a real sparse solver, 
   in Section~\ref{random results} we present timing
   results for large random (dense) matrices that are constructed to avoid the
   need for pivoting. This allows us to simulate the performance overheads of
   various pivoting techniques in the best-case where no pivot candidates are
   rejected. However, we note that when a significant number of pivot candidates
   are rejected we would expect the advantage from using compressed pivoting
   techniques to increase as threshold partial pivoting would require
   additional communication to retest these pivots later in the factorization.
   To explore the reliability of
   the compressed and restricted pivoting algorithms, in Section~\ref{real-world results} we
   present numerical stability results for a set of sparse
   problems arising from practical applications.
   
   All experiments are performed on the machine summarised in Table~\ref{irons}.

   \begin{table}[htbp]
      \caption{ \label{irons}
         Description of machine used for numerical experiments
      }
      \centering
      \begin{tabular}{ll}
          \\
         \hline
         \textbf{Processor} & 2 $\times$ Intel Xeon E5-2687W\\
         \textbf{Physical Cores} & 16 \\
         \textbf{Memory} & 64GB \\
         \textbf{Compiler} & ifort 12.1.0 \\
         \textbf{BLAS} & MKL 10.3 update 6 \\
         \textbf{L1/L2 cache (per core)} & 32KB / 256KB \\
         \textbf{L3 cache (shared)} & 20MB \\
         \textbf{Compiler flags} & \texttt{ifort -O3 -xAVX -no-prec-div -ip} \\
         \hline

      \end{tabular}
   \end{table}

   \subsection{Performance experiments on large random matrices}\label{random results}

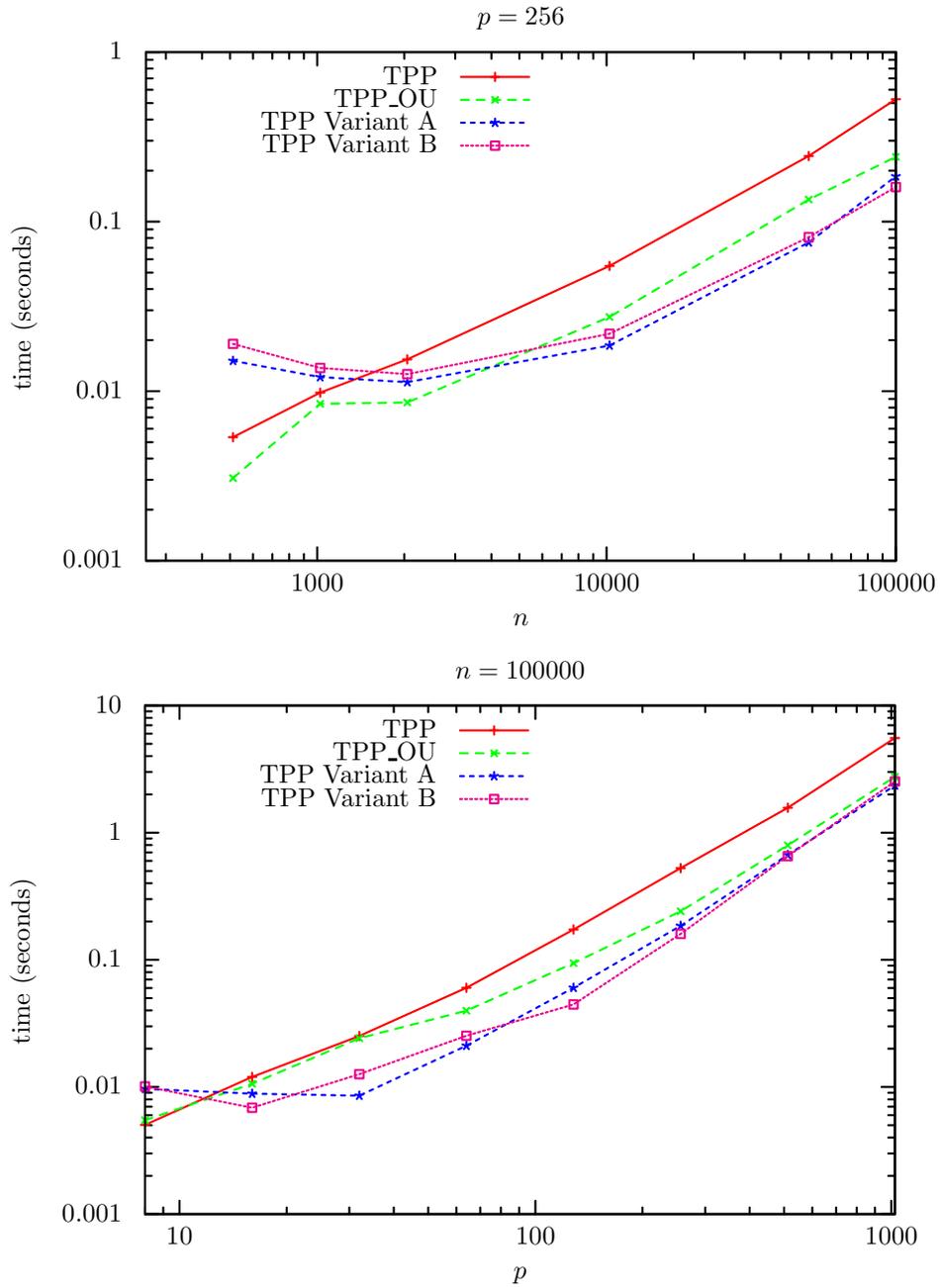
\begin{figure}
      \caption{ \label{dense tpp variants}
         Performance of variants of threshold partial pivoting using
         16 threads.
      }
      \centering
      \begin{tikzpicture}[gnuplot]
%% generated with GNUPLOT 4.6p0 (Lua 5.1; terminal rev. 99, script rev. 100)
%% Thu 07 Mar 2013 12:00:40 GMT
\path (0.000,0.000) rectangle (12.500,8.750);
\gpcolor{color=gp lt color border}
\gpsetlinetype{gp lt border}
\gpsetlinewidth{2.00}
\draw[gp path] (1.872,0.985)--(2.052,0.985);
\draw[gp path] (11.947,0.985)--(11.767,0.985);
\node[gp node right] at (1.688,0.985) { 0.001};
\draw[gp path] (1.872,1.671)--(1.962,1.671);
\draw[gp path] (11.947,1.671)--(11.857,1.671);
\draw[gp path] (1.872,2.073)--(1.962,2.073);
\draw[gp path] (11.947,2.073)--(11.857,2.073);
\draw[gp path] (1.872,2.358)--(1.962,2.358);
\draw[gp path] (11.947,2.358)--(11.857,2.358);
\draw[gp path] (1.872,2.579)--(1.962,2.579);
\draw[gp path] (11.947,2.579)--(11.857,2.579);
\draw[gp path] (1.872,2.759)--(1.962,2.759);
\draw[gp path] (11.947,2.759)--(11.857,2.759);
\draw[gp path] (1.872,2.912)--(1.962,2.912);
\draw[gp path] (11.947,2.912)--(11.857,2.912);
\draw[gp path] (1.872,3.044)--(1.962,3.044);
\draw[gp path] (11.947,3.044)--(11.857,3.044);
\draw[gp path] (1.872,3.161)--(1.962,3.161);
\draw[gp path] (11.947,3.161)--(11.857,3.161);
\draw[gp path] (1.872,3.265)--(2.052,3.265);
\draw[gp path] (11.947,3.265)--(11.767,3.265);
\node[gp node right] at (1.688,3.265) { 0.01};
\draw[gp path] (1.872,3.951)--(1.962,3.951);
\draw[gp path] (11.947,3.951)--(11.857,3.951);
\draw[gp path] (1.872,4.353)--(1.962,4.353);
\draw[gp path] (11.947,4.353)--(11.857,4.353);
\draw[gp path] (1.872,4.638)--(1.962,4.638);
\draw[gp path] (11.947,4.638)--(11.857,4.638);
\draw[gp path] (1.872,4.859)--(1.962,4.859);
\draw[gp path] (11.947,4.859)--(11.857,4.859);
\draw[gp path] (1.872,5.039)--(1.962,5.039);
\draw[gp path] (11.947,5.039)--(11.857,5.039);
\draw[gp path] (1.872,5.192)--(1.962,5.192);
\draw[gp path] (11.947,5.192)--(11.857,5.192);
\draw[gp path] (1.872,5.324)--(1.962,5.324);
\draw[gp path] (11.947,5.324)--(11.857,5.324);
\draw[gp path] (1.872,5.441)--(1.962,5.441);
\draw[gp path] (11.947,5.441)--(11.857,5.441);
\draw[gp path] (1.872,5.545)--(2.052,5.545);
\draw[gp path] (11.947,5.545)--(11.767,5.545);
\node[gp node right] at (1.688,5.545) { 0.1};
\draw[gp path] (1.872,6.231)--(1.962,6.231);
\draw[gp path] (11.947,6.231)--(11.857,6.231);
\draw[gp path] (1.872,6.633)--(1.962,6.633);
\draw[gp path] (11.947,6.633)--(11.857,6.633);
\draw[gp path] (1.872,6.918)--(1.962,6.918);
\draw[gp path] (11.947,6.918)--(11.857,6.918);
\draw[gp path] (1.872,7.139)--(1.962,7.139);
\draw[gp path] (11.947,7.139)--(11.857,7.139);
\draw[gp path] (1.872,7.319)--(1.962,7.319);
\draw[gp path] (11.947,7.319)--(11.857,7.319);
\draw[gp path] (1.872,7.472)--(1.962,7.472);
\draw[gp path] (11.947,7.472)--(11.857,7.472);
\draw[gp path] (1.872,7.604)--(1.962,7.604);
\draw[gp path] (11.947,7.604)--(11.857,7.604);
\draw[gp path] (1.872,7.721)--(1.962,7.721);
\draw[gp path] (11.947,7.721)--(11.857,7.721);
\draw[gp path] (1.872,7.825)--(2.052,7.825);
\draw[gp path] (11.947,7.825)--(11.767,7.825);
\node[gp node right] at (1.688,7.825) { 1};
\draw[gp path] (2.140,0.985)--(2.140,1.075);
\draw[gp path] (2.140,7.825)--(2.140,7.735);
\draw[gp path] (2.625,0.985)--(2.625,1.075);
\draw[gp path] (2.625,7.825)--(2.625,7.735);
\draw[gp path] (3.002,0.985)--(3.002,1.075);
\draw[gp path] (3.002,7.825)--(3.002,7.735);
\draw[gp path] (3.310,0.985)--(3.310,1.075);
\draw[gp path] (3.310,7.825)--(3.310,7.735);
\draw[gp path] (3.570,0.985)--(3.570,1.075);
\draw[gp path] (3.570,7.825)--(3.570,7.735);
\draw[gp path] (3.796,0.985)--(3.796,1.075);
\draw[gp path] (3.796,7.825)--(3.796,7.735);
\draw[gp path] (3.994,0.985)--(3.994,1.075);
\draw[gp path] (3.994,7.825)--(3.994,7.735);
\draw[gp path] (4.172,0.985)--(4.172,1.165);
\draw[gp path] (4.172,7.825)--(4.172,7.645);
\node[gp node center] at (4.172,0.677) { 1000};
\draw[gp path] (5.343,0.985)--(5.343,1.075);
\draw[gp path] (5.343,7.825)--(5.343,7.735);
\draw[gp path] (6.027,0.985)--(6.027,1.075);
\draw[gp path] (6.027,7.825)--(6.027,7.735);
\draw[gp path] (6.513,0.985)--(6.513,1.075);
\draw[gp path] (6.513,7.825)--(6.513,7.735);
\draw[gp path] (6.889,0.985)--(6.889,1.075);
\draw[gp path] (6.889,7.825)--(6.889,7.735);
\draw[gp path] (7.197,0.985)--(7.197,1.075);
\draw[gp path] (7.197,7.825)--(7.197,7.735);
\draw[gp path] (7.458,0.985)--(7.458,1.075);
\draw[gp path] (7.458,7.825)--(7.458,7.735);
\draw[gp path] (7.683,0.985)--(7.683,1.075);
\draw[gp path] (7.683,7.825)--(7.683,7.735);
\draw[gp path] (7.882,0.985)--(7.882,1.075);
\draw[gp path] (7.882,7.825)--(7.882,7.735);
\draw[gp path] (8.060,0.985)--(8.060,1.165);
\draw[gp path] (8.060,7.825)--(8.060,7.645);
\node[gp node center] at (8.060,0.677) { 10000};
\draw[gp path] (9.230,0.985)--(9.230,1.075);
\draw[gp path] (9.230,7.825)--(9.230,7.735);
\draw[gp path] (9.914,0.985)--(9.914,1.075);
\draw[gp path] (9.914,7.825)--(9.914,7.735);
\draw[gp path] (10.400,0.985)--(10.400,1.075);
\draw[gp path] (10.400,7.825)--(10.400,7.735);
\draw[gp path] (10.777,0.985)--(10.777,1.075);
\draw[gp path] (10.777,7.825)--(10.777,7.735);
\draw[gp path] (11.085,0.985)--(11.085,1.075);
\draw[gp path] (11.085,7.825)--(11.085,7.735);
\draw[gp path] (11.345,0.985)--(11.345,1.075);
\draw[gp path] (11.345,7.825)--(11.345,7.735);
\draw[gp path] (11.570,0.985)--(11.570,1.075);
\draw[gp path] (11.570,7.825)--(11.570,7.735);
\draw[gp path] (11.769,0.985)--(11.769,1.075);
\draw[gp path] (11.769,7.825)--(11.769,7.735);
\draw[gp path] (11.947,0.985)--(11.947,1.165);
\draw[gp path] (11.947,7.825)--(11.947,7.645);
\node[gp node center] at (11.947,0.677) { 100000};
\draw[gp path] (1.872,7.825)--(1.872,0.985)--(11.947,0.985)--(11.947,7.825)--cycle;
\node[gp node center,rotate=-270] at (0.246,4.405) {time (seconds)};
\node[gp node center] at (6.909,0.215) {$n$};
\node[gp node center] at (6.909,8.287) {$p=256$};
\node[gp node right] at (5.920,7.491) {TPP};
\gpcolor{color=gp lt color 0}
\gpsetlinetype{gp lt plot 0}
\draw[gp path] (6.104,7.491)--(7.020,7.491);
\draw[gp path] (3.042,2.646)--(4.212,3.245)--(5.383,3.693)--(8.100,4.949)--(10.777,6.428)%
  --(11.947,7.191);
\gpsetpointsize{4.00}
\gppoint{gp mark 1}{(3.042,2.646)}
\gppoint{gp mark 1}{(4.212,3.245)}
\gppoint{gp mark 1}{(5.383,3.693)}
\gppoint{gp mark 1}{(8.100,4.949)}
\gppoint{gp mark 1}{(10.777,6.428)}
\gppoint{gp mark 1}{(11.947,7.191)}
\gppoint{gp mark 1}{(6.562,7.491)}
\gpcolor{color=gp lt color border}
\node[gp node right] at (5.920,7.183) {TPP\_OU};
\gpcolor{color=gp lt color 1}
\gpsetlinetype{gp lt plot 1}
\draw[gp path] (6.104,7.183)--(7.020,7.183);
\draw[gp path] (3.042,2.096)--(4.212,3.095)--(5.383,3.113)--(8.100,4.263)--(10.777,5.842)%
  --(11.947,6.416);
\gppoint{gp mark 2}{(3.042,2.096)}
\gppoint{gp mark 2}{(4.212,3.095)}
\gppoint{gp mark 2}{(5.383,3.113)}
\gppoint{gp mark 2}{(8.100,4.263)}
\gppoint{gp mark 2}{(10.777,5.842)}
\gppoint{gp mark 2}{(11.947,6.416)}
\gppoint{gp mark 2}{(6.562,7.183)}
\gpcolor{color=gp lt color border}
\node[gp node right] at (5.920,6.875) {TPP Variant A};
\gpcolor{color=gp lt color 2}
\gpsetlinetype{gp lt plot 2}
\draw[gp path] (6.104,6.875)--(7.020,6.875);
\draw[gp path] (3.042,3.673)--(4.212,3.454)--(5.383,3.386)--(8.100,3.879)--(10.777,5.265)%
  --(11.947,6.154);
\gppoint{gp mark 3}{(3.042,3.673)}
\gppoint{gp mark 3}{(4.212,3.454)}
\gppoint{gp mark 3}{(5.383,3.386)}
\gppoint{gp mark 3}{(8.100,3.879)}
\gppoint{gp mark 3}{(10.777,5.265)}
\gppoint{gp mark 3}{(11.947,6.154)}
\gppoint{gp mark 3}{(6.562,6.875)}
\gpcolor{color=gp lt color border}
\node[gp node right] at (5.920,6.567) {TPP Variant B};
\gpcolor{color=gp lt color 3}
\gpsetlinetype{gp lt plot 3}
\draw[gp path] (6.104,6.567)--(7.020,6.567);
\draw[gp path] (3.042,3.901)--(4.212,3.577)--(5.383,3.494)--(8.100,4.037)--(10.777,5.336)%
  --(11.947,6.010);
\gppoint{gp mark 4}{(3.042,3.901)}
\gppoint{gp mark 4}{(4.212,3.577)}
\gppoint{gp mark 4}{(5.383,3.494)}
\gppoint{gp mark 4}{(8.100,4.037)}
\gppoint{gp mark 4}{(10.777,5.336)}
\gppoint{gp mark 4}{(11.947,6.010)}
\gppoint{gp mark 4}{(6.562,6.567)}
\gpcolor{color=gp lt color border}
\gpsetlinetype{gp lt border}
\draw[gp path] (1.872,7.825)--(1.872,0.985)--(11.947,0.985)--(11.947,7.825)--cycle;
%% coordinates of the plot area
\gpdefrectangularnode{gp plot 1}{\pgfpoint{1.872cm}{0.985cm}}{\pgfpoint{11.947cm}{7.825cm}}
\end{tikzpicture}
%% gnuplot variables
      \begin{tikzpicture}[gnuplot]
%% generated with GNUPLOT 4.6p0 (Lua 5.1; terminal rev. 99, script rev. 100)
%% Thu 07 Mar 2013 11:59:07 GMT
\path (0.000,0.000) rectangle (12.500,8.750);
\gpcolor{color=gp lt color border}
\gpsetlinetype{gp lt border}
\gpsetlinewidth{2.00}
\draw[gp path] (1.872,0.985)--(2.052,0.985);
\draw[gp path] (11.947,0.985)--(11.767,0.985);
\node[gp node right] at (1.688,0.985) { 0.001};
\draw[gp path] (1.872,1.500)--(1.962,1.500);
\draw[gp path] (11.947,1.500)--(11.857,1.500);
\draw[gp path] (1.872,1.801)--(1.962,1.801);
\draw[gp path] (11.947,1.801)--(11.857,1.801);
\draw[gp path] (1.872,2.015)--(1.962,2.015);
\draw[gp path] (11.947,2.015)--(11.857,2.015);
\draw[gp path] (1.872,2.180)--(1.962,2.180);
\draw[gp path] (11.947,2.180)--(11.857,2.180);
\draw[gp path] (1.872,2.316)--(1.962,2.316);
\draw[gp path] (11.947,2.316)--(11.857,2.316);
\draw[gp path] (1.872,2.430)--(1.962,2.430);
\draw[gp path] (11.947,2.430)--(11.857,2.430);
\draw[gp path] (1.872,2.529)--(1.962,2.529);
\draw[gp path] (11.947,2.529)--(11.857,2.529);
\draw[gp path] (1.872,2.617)--(1.962,2.617);
\draw[gp path] (11.947,2.617)--(11.857,2.617);
\draw[gp path] (1.872,2.695)--(2.052,2.695);
\draw[gp path] (11.947,2.695)--(11.767,2.695);
\node[gp node right] at (1.688,2.695) { 0.01};
\draw[gp path] (1.872,3.210)--(1.962,3.210);
\draw[gp path] (11.947,3.210)--(11.857,3.210);
\draw[gp path] (1.872,3.511)--(1.962,3.511);
\draw[gp path] (11.947,3.511)--(11.857,3.511);
\draw[gp path] (1.872,3.725)--(1.962,3.725);
\draw[gp path] (11.947,3.725)--(11.857,3.725);
\draw[gp path] (1.872,3.890)--(1.962,3.890);
\draw[gp path] (11.947,3.890)--(11.857,3.890);
\draw[gp path] (1.872,4.026)--(1.962,4.026);
\draw[gp path] (11.947,4.026)--(11.857,4.026);
\draw[gp path] (1.872,4.140)--(1.962,4.140);
\draw[gp path] (11.947,4.140)--(11.857,4.140);
\draw[gp path] (1.872,4.239)--(1.962,4.239);
\draw[gp path] (11.947,4.239)--(11.857,4.239);
\draw[gp path] (1.872,4.327)--(1.962,4.327);
\draw[gp path] (11.947,4.327)--(11.857,4.327);
\draw[gp path] (1.872,4.405)--(2.052,4.405);
\draw[gp path] (11.947,4.405)--(11.767,4.405);
\node[gp node right] at (1.688,4.405) { 0.1};
\draw[gp path] (1.872,4.920)--(1.962,4.920);
\draw[gp path] (11.947,4.920)--(11.857,4.920);
\draw[gp path] (1.872,5.221)--(1.962,5.221);
\draw[gp path] (11.947,5.221)--(11.857,5.221);
\draw[gp path] (1.872,5.435)--(1.962,5.435);
\draw[gp path] (11.947,5.435)--(11.857,5.435);
\draw[gp path] (1.872,5.600)--(1.962,5.600);
\draw[gp path] (11.947,5.600)--(11.857,5.600);
\draw[gp path] (1.872,5.736)--(1.962,5.736);
\draw[gp path] (11.947,5.736)--(11.857,5.736);
\draw[gp path] (1.872,5.850)--(1.962,5.850);
\draw[gp path] (11.947,5.850)--(11.857,5.850);
\draw[gp path] (1.872,5.949)--(1.962,5.949);
\draw[gp path] (11.947,5.949)--(11.857,5.949);
\draw[gp path] (1.872,6.037)--(1.962,6.037);
\draw[gp path] (11.947,6.037)--(11.857,6.037);
\draw[gp path] (1.872,6.115)--(2.052,6.115);
\draw[gp path] (11.947,6.115)--(11.767,6.115);
\node[gp node right] at (1.688,6.115) { 1};
\draw[gp path] (1.872,6.630)--(1.962,6.630);
\draw[gp path] (11.947,6.630)--(11.857,6.630);
\draw[gp path] (1.872,6.931)--(1.962,6.931);
\draw[gp path] (11.947,6.931)--(11.857,6.931);
\draw[gp path] (1.872,7.145)--(1.962,7.145);
\draw[gp path] (11.947,7.145)--(11.857,7.145);
\draw[gp path] (1.872,7.310)--(1.962,7.310);
\draw[gp path] (11.947,7.310)--(11.857,7.310);
\draw[gp path] (1.872,7.446)--(1.962,7.446);
\draw[gp path] (11.947,7.446)--(11.857,7.446);
\draw[gp path] (1.872,7.560)--(1.962,7.560);
\draw[gp path] (11.947,7.560)--(11.857,7.560);
\draw[gp path] (1.872,7.659)--(1.962,7.659);
\draw[gp path] (11.947,7.659)--(11.857,7.659);
\draw[gp path] (1.872,7.747)--(1.962,7.747);
\draw[gp path] (11.947,7.747)--(11.857,7.747);
\draw[gp path] (1.872,7.825)--(2.052,7.825);
\draw[gp path] (11.947,7.825)--(11.767,7.825);
\node[gp node right] at (1.688,7.825) { 10};
\draw[gp path] (1.872,0.985)--(1.872,1.075);
\draw[gp path] (1.872,7.825)--(1.872,7.735);
\draw[gp path] (2.117,0.985)--(2.117,1.075);
\draw[gp path] (2.117,7.825)--(2.117,7.735);
\draw[gp path] (2.335,0.985)--(2.335,1.165);
\draw[gp path] (2.335,7.825)--(2.335,7.645);
\node[gp node center] at (2.335,0.677) { 10};
\draw[gp path] (3.775,0.985)--(3.775,1.075);
\draw[gp path] (3.775,7.825)--(3.775,7.735);
\draw[gp path] (4.617,0.985)--(4.617,1.075);
\draw[gp path] (4.617,7.825)--(4.617,7.735);
\draw[gp path] (5.214,0.985)--(5.214,1.075);
\draw[gp path] (5.214,7.825)--(5.214,7.735);
\draw[gp path] (5.677,0.985)--(5.677,1.075);
\draw[gp path] (5.677,7.825)--(5.677,7.735);
\draw[gp path] (6.056,0.985)--(6.056,1.075);
\draw[gp path] (6.056,7.825)--(6.056,7.735);
\draw[gp path] (6.376,0.985)--(6.376,1.075);
\draw[gp path] (6.376,7.825)--(6.376,7.735);
\draw[gp path] (6.653,0.985)--(6.653,1.075);
\draw[gp path] (6.653,7.825)--(6.653,7.735);
\draw[gp path] (6.898,0.985)--(6.898,1.075);
\draw[gp path] (6.898,7.825)--(6.898,7.735);
\draw[gp path] (7.117,0.985)--(7.117,1.165);
\draw[gp path] (7.117,7.825)--(7.117,7.645);
\node[gp node center] at (7.117,0.677) { 100};
\draw[gp path] (8.556,0.985)--(8.556,1.075);
\draw[gp path] (8.556,7.825)--(8.556,7.735);
\draw[gp path] (9.398,0.985)--(9.398,1.075);
\draw[gp path] (9.398,7.825)--(9.398,7.735);
\draw[gp path] (9.995,0.985)--(9.995,1.075);
\draw[gp path] (9.995,7.825)--(9.995,7.735);
\draw[gp path] (10.458,0.985)--(10.458,1.075);
\draw[gp path] (10.458,7.825)--(10.458,7.735);
\draw[gp path] (10.837,0.985)--(10.837,1.075);
\draw[gp path] (10.837,7.825)--(10.837,7.735);
\draw[gp path] (11.157,0.985)--(11.157,1.075);
\draw[gp path] (11.157,7.825)--(11.157,7.735);
\draw[gp path] (11.434,0.985)--(11.434,1.075);
\draw[gp path] (11.434,7.825)--(11.434,7.735);
\draw[gp path] (11.679,0.985)--(11.679,1.075);
\draw[gp path] (11.679,7.825)--(11.679,7.735);
\draw[gp path] (11.898,0.985)--(11.898,1.165);
\draw[gp path] (11.898,7.825)--(11.898,7.645);
\node[gp node center] at (11.898,0.677) { 1000};
\draw[gp path] (1.872,7.825)--(1.872,0.985)--(11.947,0.985)--(11.947,7.825)--cycle;
\node[gp node center,rotate=-270] at (0.246,4.405) {time (seconds)};
\node[gp node center] at (6.909,0.215) {$p$};
\node[gp node center] at (6.909,8.287) {$n=100000$};
\node[gp node right] at (5.920,7.491) {TPP};
\gpcolor{color=gp lt color 0}
\gpsetlinetype{gp lt plot 0}
\draw[gp path] (6.104,7.491)--(7.020,7.491);
\draw[gp path] (1.872,2.185)--(3.311,2.830)--(4.751,3.381)--(6.190,4.029)--(7.629,4.812)%
  --(9.068,5.639)--(10.508,6.450)--(11.947,7.388);
\gpsetpointsize{4.00}
\gppoint{gp mark 1}{(1.872,2.185)}
\gppoint{gp mark 1}{(3.311,2.830)}
\gppoint{gp mark 1}{(4.751,3.381)}
\gppoint{gp mark 1}{(6.190,4.029)}
\gppoint{gp mark 1}{(7.629,4.812)}
\gppoint{gp mark 1}{(9.068,5.639)}
\gppoint{gp mark 1}{(10.508,6.450)}
\gppoint{gp mark 1}{(11.947,7.388)}
\gppoint{gp mark 1}{(6.562,7.491)}
\gpcolor{color=gp lt color border}
\node[gp node right] at (5.920,7.183) {TPP\_OU};
\gpcolor{color=gp lt color 1}
\gpsetlinetype{gp lt plot 1}
\draw[gp path] (6.104,7.183)--(7.020,7.183);
\draw[gp path] (1.872,2.248)--(3.311,2.738)--(4.751,3.351)--(6.190,3.721)--(7.629,4.360)%
  --(9.068,5.058)--(10.508,5.946)--(11.947,6.866);
\gppoint{gp mark 2}{(1.872,2.248)}
\gppoint{gp mark 2}{(3.311,2.738)}
\gppoint{gp mark 2}{(4.751,3.351)}
\gppoint{gp mark 2}{(6.190,3.721)}
\gppoint{gp mark 2}{(7.629,4.360)}
\gppoint{gp mark 2}{(9.068,5.058)}
\gppoint{gp mark 2}{(10.508,5.946)}
\gppoint{gp mark 2}{(11.947,6.866)}
\gppoint{gp mark 2}{(6.562,7.183)}
\gpcolor{color=gp lt color border}
\node[gp node right] at (5.920,6.875) {TPP Variant A};
\gpcolor{color=gp lt color 2}
\gpsetlinetype{gp lt plot 2}
\draw[gp path] (6.104,6.875)--(7.020,6.875);
\draw[gp path] (1.872,2.670)--(3.311,2.605)--(4.751,2.578)--(6.190,3.250)--(7.629,4.029)%
  --(9.068,4.862)--(10.508,5.816)--(11.947,6.750);
\gppoint{gp mark 3}{(1.872,2.670)}
\gppoint{gp mark 3}{(3.311,2.605)}
\gppoint{gp mark 3}{(4.751,2.578)}
\gppoint{gp mark 3}{(6.190,3.250)}
\gppoint{gp mark 3}{(7.629,4.029)}
\gppoint{gp mark 3}{(9.068,4.862)}
\gppoint{gp mark 3}{(10.508,5.816)}
\gppoint{gp mark 3}{(11.947,6.750)}
\gppoint{gp mark 3}{(6.562,6.875)}
\gpcolor{color=gp lt color border}
\node[gp node right] at (5.920,6.567) {TPP Variant B};
\gpcolor{color=gp lt color 3}
\gpsetlinetype{gp lt plot 3}
\draw[gp path] (6.104,6.567)--(7.020,6.567);
\draw[gp path] (1.872,2.702)--(3.311,2.416)--(4.751,2.867)--(6.190,3.384)--(7.629,3.805)%
  --(9.068,4.754)--(10.508,5.800)--(11.947,6.804);
\gppoint{gp mark 4}{(1.872,2.702)}
\gppoint{gp mark 4}{(3.311,2.416)}
\gppoint{gp mark 4}{(4.751,2.867)}
\gppoint{gp mark 4}{(6.190,3.384)}
\gppoint{gp mark 4}{(7.629,3.805)}
\gppoint{gp mark 4}{(9.068,4.754)}
\gppoint{gp mark 4}{(10.508,5.800)}
\gppoint{gp mark 4}{(11.947,6.804)}
\gppoint{gp mark 4}{(6.562,6.567)}
\gpcolor{color=gp lt color border}
\gpsetlinetype{gp lt border}
\draw[gp path] (1.872,7.825)--(1.872,0.985)--(11.947,0.985)--(11.947,7.825)--cycle;
%% coordinates of the plot area
\gpdefrectangularnode{gp plot 1}{\pgfpoint{1.872cm}{0.985cm}}{\pgfpoint{11.947cm}{7.825cm}}
\end{tikzpicture}
%% gnuplot variables
   \end{figure}

   Figure~\ref{dense tpp variants} shows the performance characteristics for
   the threshold partial pivoting variants 
   discussed in Section~\ref{TPP:variants}. For the TPP\_OU
   variant we use block column size $nbi=16$. The top graph shows a slice through the $(n,p)$
   parameter space for large fixed $p$, while the bottom graph shows a slice
   for large fixed $n$.

   As we might expect, the figure shows that for small $n$ the TPP\_OU
   approach is the fastest because it avoids the communication overheads
   of the inner loops
   inherent in Variants A and B. For larger $n$, there is sufficient work to
   amortize such overheads and Variants A and B perform best, with a slight
   performance advantage for Variant B. As
   $p$ increases, more time is spent in the outer update, and there is therefore
   little to choose between the parallel implementations.

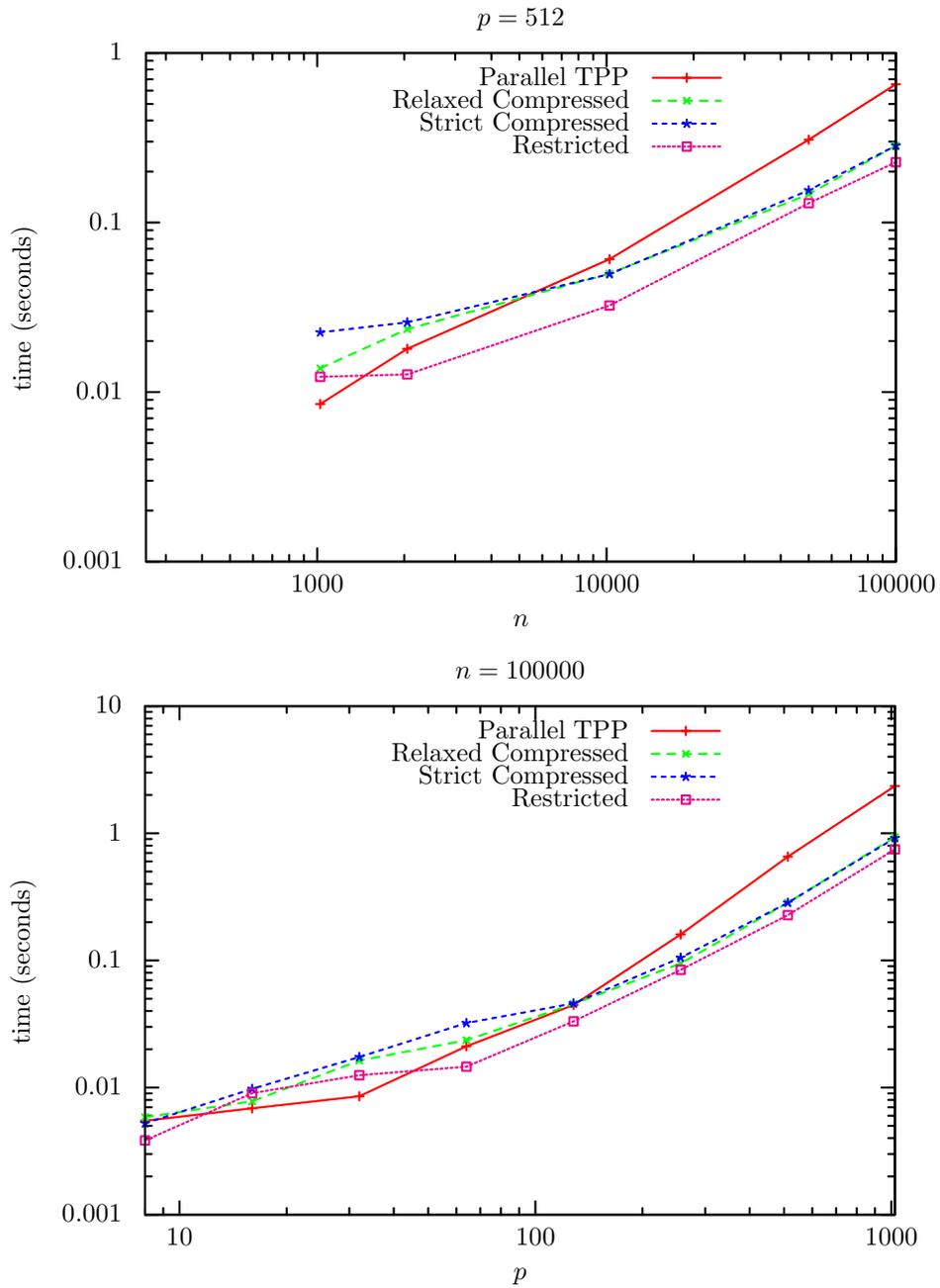
\begin{figure}
      \caption{ \label{dense parallel}
         Performance of parallel pivoting schemes using 16 threads. For
         TPP the best variant
         is used for each combination of $n$ and $p$.
      }
      \centering
      \begin{tikzpicture}[gnuplot]
%% generated with GNUPLOT 4.6p0 (Lua 5.1; terminal rev. 99, script rev. 100)
%% Thu 07 Mar 2013 12:06:38 GMT
\path (0.000,0.000) rectangle (12.500,8.750);
\gpcolor{color=gp lt color border}
\gpsetlinetype{gp lt border}
\gpsetlinewidth{2.00}
\draw[gp path] (1.872,0.985)--(2.052,0.985);
\draw[gp path] (11.947,0.985)--(11.767,0.985);
\node[gp node right] at (1.688,0.985) { 0.001};
\draw[gp path] (1.872,1.671)--(1.962,1.671);
\draw[gp path] (11.947,1.671)--(11.857,1.671);
\draw[gp path] (1.872,2.073)--(1.962,2.073);
\draw[gp path] (11.947,2.073)--(11.857,2.073);
\draw[gp path] (1.872,2.358)--(1.962,2.358);
\draw[gp path] (11.947,2.358)--(11.857,2.358);
\draw[gp path] (1.872,2.579)--(1.962,2.579);
\draw[gp path] (11.947,2.579)--(11.857,2.579);
\draw[gp path] (1.872,2.759)--(1.962,2.759);
\draw[gp path] (11.947,2.759)--(11.857,2.759);
\draw[gp path] (1.872,2.912)--(1.962,2.912);
\draw[gp path] (11.947,2.912)--(11.857,2.912);
\draw[gp path] (1.872,3.044)--(1.962,3.044);
\draw[gp path] (11.947,3.044)--(11.857,3.044);
\draw[gp path] (1.872,3.161)--(1.962,3.161);
\draw[gp path] (11.947,3.161)--(11.857,3.161);
\draw[gp path] (1.872,3.265)--(2.052,3.265);
\draw[gp path] (11.947,3.265)--(11.767,3.265);
\node[gp node right] at (1.688,3.265) { 0.01};
\draw[gp path] (1.872,3.951)--(1.962,3.951);
\draw[gp path] (11.947,3.951)--(11.857,3.951);
\draw[gp path] (1.872,4.353)--(1.962,4.353);
\draw[gp path] (11.947,4.353)--(11.857,4.353);
\draw[gp path] (1.872,4.638)--(1.962,4.638);
\draw[gp path] (11.947,4.638)--(11.857,4.638);
\draw[gp path] (1.872,4.859)--(1.962,4.859);
\draw[gp path] (11.947,4.859)--(11.857,4.859);
\draw[gp path] (1.872,5.039)--(1.962,5.039);
\draw[gp path] (11.947,5.039)--(11.857,5.039);
\draw[gp path] (1.872,5.192)--(1.962,5.192);
\draw[gp path] (11.947,5.192)--(11.857,5.192);
\draw[gp path] (1.872,5.324)--(1.962,5.324);
\draw[gp path] (11.947,5.324)--(11.857,5.324);
\draw[gp path] (1.872,5.441)--(1.962,5.441);
\draw[gp path] (11.947,5.441)--(11.857,5.441);
\draw[gp path] (1.872,5.545)--(2.052,5.545);
\draw[gp path] (11.947,5.545)--(11.767,5.545);
\node[gp node right] at (1.688,5.545) { 0.1};
\draw[gp path] (1.872,6.231)--(1.962,6.231);
\draw[gp path] (11.947,6.231)--(11.857,6.231);
\draw[gp path] (1.872,6.633)--(1.962,6.633);
\draw[gp path] (11.947,6.633)--(11.857,6.633);
\draw[gp path] (1.872,6.918)--(1.962,6.918);
\draw[gp path] (11.947,6.918)--(11.857,6.918);
\draw[gp path] (1.872,7.139)--(1.962,7.139);
\draw[gp path] (11.947,7.139)--(11.857,7.139);
\draw[gp path] (1.872,7.319)--(1.962,7.319);
\draw[gp path] (11.947,7.319)--(11.857,7.319);
\draw[gp path] (1.872,7.472)--(1.962,7.472);
\draw[gp path] (11.947,7.472)--(11.857,7.472);
\draw[gp path] (1.872,7.604)--(1.962,7.604);
\draw[gp path] (11.947,7.604)--(11.857,7.604);
\draw[gp path] (1.872,7.721)--(1.962,7.721);
\draw[gp path] (11.947,7.721)--(11.857,7.721);
\draw[gp path] (1.872,7.825)--(2.052,7.825);
\draw[gp path] (11.947,7.825)--(11.767,7.825);
\node[gp node right] at (1.688,7.825) { 1};
\draw[gp path] (2.140,0.985)--(2.140,1.075);
\draw[gp path] (2.140,7.825)--(2.140,7.735);
\draw[gp path] (2.625,0.985)--(2.625,1.075);
\draw[gp path] (2.625,7.825)--(2.625,7.735);
\draw[gp path] (3.002,0.985)--(3.002,1.075);
\draw[gp path] (3.002,7.825)--(3.002,7.735);
\draw[gp path] (3.310,0.985)--(3.310,1.075);
\draw[gp path] (3.310,7.825)--(3.310,7.735);
\draw[gp path] (3.570,0.985)--(3.570,1.075);
\draw[gp path] (3.570,7.825)--(3.570,7.735);
\draw[gp path] (3.796,0.985)--(3.796,1.075);
\draw[gp path] (3.796,7.825)--(3.796,7.735);
\draw[gp path] (3.994,0.985)--(3.994,1.075);
\draw[gp path] (3.994,7.825)--(3.994,7.735);
\draw[gp path] (4.172,0.985)--(4.172,1.165);
\draw[gp path] (4.172,7.825)--(4.172,7.645);
\node[gp node center] at (4.172,0.677) { 1000};
\draw[gp path] (5.343,0.985)--(5.343,1.075);
\draw[gp path] (5.343,7.825)--(5.343,7.735);
\draw[gp path] (6.027,0.985)--(6.027,1.075);
\draw[gp path] (6.027,7.825)--(6.027,7.735);
\draw[gp path] (6.513,0.985)--(6.513,1.075);
\draw[gp path] (6.513,7.825)--(6.513,7.735);
\draw[gp path] (6.889,0.985)--(6.889,1.075);
\draw[gp path] (6.889,7.825)--(6.889,7.735);
\draw[gp path] (7.197,0.985)--(7.197,1.075);
\draw[gp path] (7.197,7.825)--(7.197,7.735);
\draw[gp path] (7.458,0.985)--(7.458,1.075);
\draw[gp path] (7.458,7.825)--(7.458,7.735);
\draw[gp path] (7.683,0.985)--(7.683,1.075);
\draw[gp path] (7.683,7.825)--(7.683,7.735);
\draw[gp path] (7.882,0.985)--(7.882,1.075);
\draw[gp path] (7.882,7.825)--(7.882,7.735);
\draw[gp path] (8.060,0.985)--(8.060,1.165);
\draw[gp path] (8.060,7.825)--(8.060,7.645);
\node[gp node center] at (8.060,0.677) { 10000};
\draw[gp path] (9.230,0.985)--(9.230,1.075);
\draw[gp path] (9.230,7.825)--(9.230,7.735);
\draw[gp path] (9.914,0.985)--(9.914,1.075);
\draw[gp path] (9.914,7.825)--(9.914,7.735);
\draw[gp path] (10.400,0.985)--(10.400,1.075);
\draw[gp path] (10.400,7.825)--(10.400,7.735);
\draw[gp path] (10.777,0.985)--(10.777,1.075);
\draw[gp path] (10.777,7.825)--(10.777,7.735);
\draw[gp path] (11.085,0.985)--(11.085,1.075);
\draw[gp path] (11.085,7.825)--(11.085,7.735);
\draw[gp path] (11.345,0.985)--(11.345,1.075);
\draw[gp path] (11.345,7.825)--(11.345,7.735);
\draw[gp path] (11.570,0.985)--(11.570,1.075);
\draw[gp path] (11.570,7.825)--(11.570,7.735);
\draw[gp path] (11.769,0.985)--(11.769,1.075);
\draw[gp path] (11.769,7.825)--(11.769,7.735);
\draw[gp path] (11.947,0.985)--(11.947,1.165);
\draw[gp path] (11.947,7.825)--(11.947,7.645);
\node[gp node center] at (11.947,0.677) { 100000};
\draw[gp path] (1.872,7.825)--(1.872,0.985)--(11.947,0.985)--(11.947,7.825)--cycle;
\node[gp node center,rotate=-270] at (0.246,4.405) {time (seconds)};
\node[gp node center] at (6.909,0.215) {$n$};
\node[gp node center] at (6.909,8.287) {$p=512$};
\node[gp node right] at (8.496,7.491) {Parallel TPP};
\gpcolor{color=gp lt color 0}
\gpsetlinetype{gp lt plot 0}
\draw[gp path] (8.680,7.491)--(9.596,7.491);
\draw[gp path] (4.212,3.103)--(5.383,3.847)--(8.100,5.052)--(10.777,6.659)--(11.947,7.405);
\gpsetpointsize{4.00}
\gppoint{gp mark 1}{(4.212,3.103)}
\gppoint{gp mark 1}{(5.383,3.847)}
\gppoint{gp mark 1}{(8.100,5.052)}
\gppoint{gp mark 1}{(10.777,6.659)}
\gppoint{gp mark 1}{(11.947,7.405)}
\gppoint{gp mark 1}{(9.138,7.491)}
\gpcolor{color=gp lt color border}
\node[gp node right] at (8.496,7.183) {Relaxed Compressed};
\gpcolor{color=gp lt color 1}
\gpsetlinetype{gp lt plot 1}
\draw[gp path] (8.680,7.183)--(9.596,7.183);
\draw[gp path] (4.212,3.584)--(5.383,4.111)--(8.100,4.857)--(10.777,5.926)--(11.947,6.579);
\gppoint{gp mark 2}{(4.212,3.584)}
\gppoint{gp mark 2}{(5.383,4.111)}
\gppoint{gp mark 2}{(8.100,4.857)}
\gppoint{gp mark 2}{(10.777,5.926)}
\gppoint{gp mark 2}{(11.947,6.579)}
\gppoint{gp mark 2}{(9.138,7.183)}
\gpcolor{color=gp lt color border}
\node[gp node right] at (8.496,6.875) {Strict Compressed};
\gpcolor{color=gp lt color 2}
\gpsetlinetype{gp lt plot 2}
\draw[gp path] (8.680,6.875)--(9.596,6.875);
\draw[gp path] (4.212,4.068)--(5.383,4.203)--(8.100,4.851)--(10.777,5.979)--(11.947,6.579);
\gppoint{gp mark 3}{(4.212,4.068)}
\gppoint{gp mark 3}{(5.383,4.203)}
\gppoint{gp mark 3}{(8.100,4.851)}
\gppoint{gp mark 3}{(10.777,5.979)}
\gppoint{gp mark 3}{(11.947,6.579)}
\gppoint{gp mark 3}{(9.138,6.875)}
\gpcolor{color=gp lt color border}
\node[gp node right] at (8.496,6.567) {Restricted};
\gpcolor{color=gp lt color 3}
\gpsetlinetype{gp lt plot 3}
\draw[gp path] (8.680,6.567)--(9.596,6.567);
\draw[gp path] (4.212,3.470)--(5.383,3.502)--(8.100,4.429)--(10.777,5.805)--(11.947,6.357);
\gppoint{gp mark 4}{(4.212,3.470)}
\gppoint{gp mark 4}{(5.383,3.502)}
\gppoint{gp mark 4}{(8.100,4.429)}
\gppoint{gp mark 4}{(10.777,5.805)}
\gppoint{gp mark 4}{(11.947,6.357)}
\gppoint{gp mark 4}{(9.138,6.567)}
\gpcolor{color=gp lt color border}
\gpsetlinetype{gp lt border}
\draw[gp path] (1.872,7.825)--(1.872,0.985)--(11.947,0.985)--(11.947,7.825)--cycle;
%% coordinates of the plot area
\gpdefrectangularnode{gp plot 1}{\pgfpoint{1.872cm}{0.985cm}}{\pgfpoint{11.947cm}{7.825cm}}
\end{tikzpicture}
%% gnuplot variables
      \begin{tikzpicture}[gnuplot]
%% generated with GNUPLOT 4.6p0 (Lua 5.1; terminal rev. 99, script rev. 100)
%% Thu 07 Mar 2013 12:10:24 GMT
\path (0.000,0.000) rectangle (12.500,8.750);
\gpcolor{color=gp lt color border}
\gpsetlinetype{gp lt border}
\gpsetlinewidth{2.00}
\draw[gp path] (1.872,0.985)--(2.052,0.985);
\draw[gp path] (11.947,0.985)--(11.767,0.985);
\node[gp node right] at (1.688,0.985) { 0.001};
\draw[gp path] (1.872,1.500)--(1.962,1.500);
\draw[gp path] (11.947,1.500)--(11.857,1.500);
\draw[gp path] (1.872,1.801)--(1.962,1.801);
\draw[gp path] (11.947,1.801)--(11.857,1.801);
\draw[gp path] (1.872,2.015)--(1.962,2.015);
\draw[gp path] (11.947,2.015)--(11.857,2.015);
\draw[gp path] (1.872,2.180)--(1.962,2.180);
\draw[gp path] (11.947,2.180)--(11.857,2.180);
\draw[gp path] (1.872,2.316)--(1.962,2.316);
\draw[gp path] (11.947,2.316)--(11.857,2.316);
\draw[gp path] (1.872,2.430)--(1.962,2.430);
\draw[gp path] (11.947,2.430)--(11.857,2.430);
\draw[gp path] (1.872,2.529)--(1.962,2.529);
\draw[gp path] (11.947,2.529)--(11.857,2.529);
\draw[gp path] (1.872,2.617)--(1.962,2.617);
\draw[gp path] (11.947,2.617)--(11.857,2.617);
\draw[gp path] (1.872,2.695)--(2.052,2.695);
\draw[gp path] (11.947,2.695)--(11.767,2.695);
\node[gp node right] at (1.688,2.695) { 0.01};
\draw[gp path] (1.872,3.210)--(1.962,3.210);
\draw[gp path] (11.947,3.210)--(11.857,3.210);
\draw[gp path] (1.872,3.511)--(1.962,3.511);
\draw[gp path] (11.947,3.511)--(11.857,3.511);
\draw[gp path] (1.872,3.725)--(1.962,3.725);
\draw[gp path] (11.947,3.725)--(11.857,3.725);
\draw[gp path] (1.872,3.890)--(1.962,3.890);
\draw[gp path] (11.947,3.890)--(11.857,3.890);
\draw[gp path] (1.872,4.026)--(1.962,4.026);
\draw[gp path] (11.947,4.026)--(11.857,4.026);
\draw[gp path] (1.872,4.140)--(1.962,4.140);
\draw[gp path] (11.947,4.140)--(11.857,4.140);
\draw[gp path] (1.872,4.239)--(1.962,4.239);
\draw[gp path] (11.947,4.239)--(11.857,4.239);
\draw[gp path] (1.872,4.327)--(1.962,4.327);
\draw[gp path] (11.947,4.327)--(11.857,4.327);
\draw[gp path] (1.872,4.405)--(2.052,4.405);
\draw[gp path] (11.947,4.405)--(11.767,4.405);
\node[gp node right] at (1.688,4.405) { 0.1};
\draw[gp path] (1.872,4.920)--(1.962,4.920);
\draw[gp path] (11.947,4.920)--(11.857,4.920);
\draw[gp path] (1.872,5.221)--(1.962,5.221);
\draw[gp path] (11.947,5.221)--(11.857,5.221);
\draw[gp path] (1.872,5.435)--(1.962,5.435);
\draw[gp path] (11.947,5.435)--(11.857,5.435);
\draw[gp path] (1.872,5.600)--(1.962,5.600);
\draw[gp path] (11.947,5.600)--(11.857,5.600);
\draw[gp path] (1.872,5.736)--(1.962,5.736);
\draw[gp path] (11.947,5.736)--(11.857,5.736);
\draw[gp path] (1.872,5.850)--(1.962,5.850);
\draw[gp path] (11.947,5.850)--(11.857,5.850);
\draw[gp path] (1.872,5.949)--(1.962,5.949);
\draw[gp path] (11.947,5.949)--(11.857,5.949);
\draw[gp path] (1.872,6.037)--(1.962,6.037);
\draw[gp path] (11.947,6.037)--(11.857,6.037);
\draw[gp path] (1.872,6.115)--(2.052,6.115);
\draw[gp path] (11.947,6.115)--(11.767,6.115);
\node[gp node right] at (1.688,6.115) { 1};
\draw[gp path] (1.872,6.630)--(1.962,6.630);
\draw[gp path] (11.947,6.630)--(11.857,6.630);
\draw[gp path] (1.872,6.931)--(1.962,6.931);
\draw[gp path] (11.947,6.931)--(11.857,6.931);
\draw[gp path] (1.872,7.145)--(1.962,7.145);
\draw[gp path] (11.947,7.145)--(11.857,7.145);
\draw[gp path] (1.872,7.310)--(1.962,7.310);
\draw[gp path] (11.947,7.310)--(11.857,7.310);
\draw[gp path] (1.872,7.446)--(1.962,7.446);
\draw[gp path] (11.947,7.446)--(11.857,7.446);
\draw[gp path] (1.872,7.560)--(1.962,7.560);
\draw[gp path] (11.947,7.560)--(11.857,7.560);
\draw[gp path] (1.872,7.659)--(1.962,7.659);
\draw[gp path] (11.947,7.659)--(11.857,7.659);
\draw[gp path] (1.872,7.747)--(1.962,7.747);
\draw[gp path] (11.947,7.747)--(11.857,7.747);
\draw[gp path] (1.872,7.825)--(2.052,7.825);
\draw[gp path] (11.947,7.825)--(11.767,7.825);
\node[gp node right] at (1.688,7.825) { 10};
\draw[gp path] (1.872,0.985)--(1.872,1.075);
\draw[gp path] (1.872,7.825)--(1.872,7.735);
\draw[gp path] (2.117,0.985)--(2.117,1.075);
\draw[gp path] (2.117,7.825)--(2.117,7.735);
\draw[gp path] (2.335,0.985)--(2.335,1.165);
\draw[gp path] (2.335,7.825)--(2.335,7.645);
\node[gp node center] at (2.335,0.677) { 10};
\draw[gp path] (3.775,0.985)--(3.775,1.075);
\draw[gp path] (3.775,7.825)--(3.775,7.735);
\draw[gp path] (4.617,0.985)--(4.617,1.075);
\draw[gp path] (4.617,7.825)--(4.617,7.735);
\draw[gp path] (5.214,0.985)--(5.214,1.075);
\draw[gp path] (5.214,7.825)--(5.214,7.735);
\draw[gp path] (5.677,0.985)--(5.677,1.075);
\draw[gp path] (5.677,7.825)--(5.677,7.735);
\draw[gp path] (6.056,0.985)--(6.056,1.075);
\draw[gp path] (6.056,7.825)--(6.056,7.735);
\draw[gp path] (6.376,0.985)--(6.376,1.075);
\draw[gp path] (6.376,7.825)--(6.376,7.735);
\draw[gp path] (6.653,0.985)--(6.653,1.075);
\draw[gp path] (6.653,7.825)--(6.653,7.735);
\draw[gp path] (6.898,0.985)--(6.898,1.075);
\draw[gp path] (6.898,7.825)--(6.898,7.735);
\draw[gp path] (7.117,0.985)--(7.117,1.165);
\draw[gp path] (7.117,7.825)--(7.117,7.645);
\node[gp node center] at (7.117,0.677) { 100};
\draw[gp path] (8.556,0.985)--(8.556,1.075);
\draw[gp path] (8.556,7.825)--(8.556,7.735);
\draw[gp path] (9.398,0.985)--(9.398,1.075);
\draw[gp path] (9.398,7.825)--(9.398,7.735);
\draw[gp path] (9.995,0.985)--(9.995,1.075);
\draw[gp path] (9.995,7.825)--(9.995,7.735);
\draw[gp path] (10.458,0.985)--(10.458,1.075);
\draw[gp path] (10.458,7.825)--(10.458,7.735);
\draw[gp path] (10.837,0.985)--(10.837,1.075);
\draw[gp path] (10.837,7.825)--(10.837,7.735);
\draw[gp path] (11.157,0.985)--(11.157,1.075);
\draw[gp path] (11.157,7.825)--(11.157,7.735);
\draw[gp path] (11.434,0.985)--(11.434,1.075);
\draw[gp path] (11.434,7.825)--(11.434,7.735);
\draw[gp path] (11.679,0.985)--(11.679,1.075);
\draw[gp path] (11.679,7.825)--(11.679,7.735);
\draw[gp path] (11.898,0.985)--(11.898,1.165);
\draw[gp path] (11.898,7.825)--(11.898,7.645);
\node[gp node center] at (11.898,0.677) { 1000};
\draw[gp path] (1.872,7.825)--(1.872,0.985)--(11.947,0.985)--(11.947,7.825)--cycle;
\node[gp node center,rotate=-270] at (0.246,4.405) {time (seconds)};
\node[gp node center] at (6.909,0.215) {$p$};
\node[gp node center] at (6.909,8.287) {$n=100000$};
\node[gp node right] at (8.496,7.491) {Parallel TPP};
\gpcolor{color=gp lt color 0}
\gpsetlinetype{gp lt plot 0}
\draw[gp path] (8.680,7.491)--(9.596,7.491);
\draw[gp path] (1.872,2.248)--(3.311,2.416)--(4.751,2.578)--(6.190,3.250)--(7.629,3.805)%
  --(9.068,4.754)--(10.508,5.800)--(11.947,6.750);
\gpsetpointsize{4.00}
\gppoint{gp mark 1}{(1.872,2.248)}
\gppoint{gp mark 1}{(3.311,2.416)}
\gppoint{gp mark 1}{(4.751,2.578)}
\gppoint{gp mark 1}{(6.190,3.250)}
\gppoint{gp mark 1}{(7.629,3.805)}
\gppoint{gp mark 1}{(9.068,4.754)}
\gppoint{gp mark 1}{(10.508,5.800)}
\gppoint{gp mark 1}{(11.947,6.750)}
\gppoint{gp mark 1}{(9.138,7.491)}
\gpcolor{color=gp lt color border}
\node[gp node right] at (8.496,7.183) {Relaxed Compressed};
\gpcolor{color=gp lt color 1}
\gpsetlinetype{gp lt plot 1}
\draw[gp path] (8.680,7.183)--(9.596,7.183);
\draw[gp path] (1.872,2.297)--(3.311,2.510)--(4.751,3.058)--(6.190,3.333)--(7.629,3.819)%
  --(9.068,4.365)--(10.508,5.180)--(11.947,6.064);
\gppoint{gp mark 2}{(1.872,2.297)}
\gppoint{gp mark 2}{(3.311,2.510)}
\gppoint{gp mark 2}{(4.751,3.058)}
\gppoint{gp mark 2}{(6.190,3.333)}
\gppoint{gp mark 2}{(7.629,3.819)}
\gppoint{gp mark 2}{(9.068,4.365)}
\gppoint{gp mark 2}{(10.508,5.180)}
\gppoint{gp mark 2}{(11.947,6.064)}
\gppoint{gp mark 2}{(9.138,7.183)}
\gpcolor{color=gp lt color border}
\node[gp node right] at (8.496,6.875) {Strict Compressed};
\gpcolor{color=gp lt color 2}
\gpsetlinetype{gp lt plot 2}
\draw[gp path] (8.680,6.875)--(9.596,6.875);
\draw[gp path] (1.872,2.214)--(3.311,2.677)--(4.751,3.106)--(6.190,3.563)--(7.629,3.827)%
  --(9.068,4.441)--(10.508,5.180)--(11.947,6.051);
\gppoint{gp mark 3}{(1.872,2.214)}
\gppoint{gp mark 3}{(3.311,2.677)}
\gppoint{gp mark 3}{(4.751,3.106)}
\gppoint{gp mark 3}{(6.190,3.563)}
\gppoint{gp mark 3}{(7.629,3.827)}
\gppoint{gp mark 3}{(9.068,4.441)}
\gppoint{gp mark 3}{(10.508,5.180)}
\gppoint{gp mark 3}{(11.947,6.051)}
\gppoint{gp mark 3}{(9.138,6.875)}
\gpcolor{color=gp lt color border}
\node[gp node right] at (8.496,6.567) {Restricted};
\gpcolor{color=gp lt color 3}
\gpsetlinetype{gp lt plot 3}
\draw[gp path] (8.680,6.567)--(9.596,6.567);
\draw[gp path] (1.872,1.984)--(3.311,2.619)--(4.751,2.861)--(6.190,2.976)--(7.629,3.586)%
  --(9.068,4.278)--(10.508,5.014)--(11.947,5.899);
\gppoint{gp mark 4}{(1.872,1.984)}
\gppoint{gp mark 4}{(3.311,2.619)}
\gppoint{gp mark 4}{(4.751,2.861)}
\gppoint{gp mark 4}{(6.190,2.976)}
\gppoint{gp mark 4}{(7.629,3.586)}
\gppoint{gp mark 4}{(9.068,4.278)}
\gppoint{gp mark 4}{(10.508,5.014)}
\gppoint{gp mark 4}{(11.947,5.899)}
\gppoint{gp mark 4}{(9.138,6.567)}
\gpcolor{color=gp lt color border}
\gpsetlinetype{gp lt border}
\draw[gp path] (1.872,7.825)--(1.872,0.985)--(11.947,0.985)--(11.947,7.825)--cycle;
%% coordinates of the plot area
\gpdefrectangularnode{gp plot 1}{\pgfpoint{1.872cm}{0.985cm}}{\pgfpoint{11.947cm}{7.825cm}}
\end{tikzpicture}
%% gnuplot variables
   \end{figure}

   Figure~\ref{dense parallel} compares the best (on a case-by-case
   basis) variant of threshold partial pivoting with the two proposed
   compressed schemes and restricted pivoting. Note that we cannot expect either
   of the compressed schemes to be faster than restricted pivoting because, as
   already noted, they  always
   perform more operations and more communication.
   It is clear that as $n$ and $p$ increase, the compressed pivoting techniques
   substantially outperform the best threshold partial pivoting variant. For
   large $n$ and $p$ they are over twice as fast.
   Further, they almost approach the performance of restricted pivoting while
   offering better numerical stability, as will be demonstrated
   in the next section.

   \subsection{Real-world numerical stability}\label{real-world results}

   We present results for two sets of 25 sparse indefinite problems drawn from the University of
   Florida Sparse Matrix Collection \cite{Davis2011}. Test Set 1 consists of a
   selection of general problems of order at least $50 000$,
   while Test Set 2 consists
   of problems where threshold partial pivoting leads to a
   significant number of delayed pivots (these problems are selected from
   those surveyed in our recent study \cite{hosc:2013} on tough
   indefinite systems). The problems are
   scaled using a weighted matching approach (as implemented by
   {\tt MC64} \cite{duko:2001}) and solved with a modified version of our
   sparse direct solver {\tt HSL\_MA97} \cite{hosc:2011a}. In all the tests,
   we use the default threshold parameter $u = 0.01$.
   By default, {\tt HSL\_MA97} chooses between using an approximate minimum degree ordering
   and a nested dissection ordering (the choice is made on the basis
   of the order of the matrix and its density). However, it also
   offers a matching-based ordering. For tough indefinite problems,
   matching-based orderings 
   can substantially reduce the number of delayed pivots albeit at the
   possible cost of additional operations and denser factors (see
   \cite{hosc:2013,scwh:2007}).

\begin{figure}
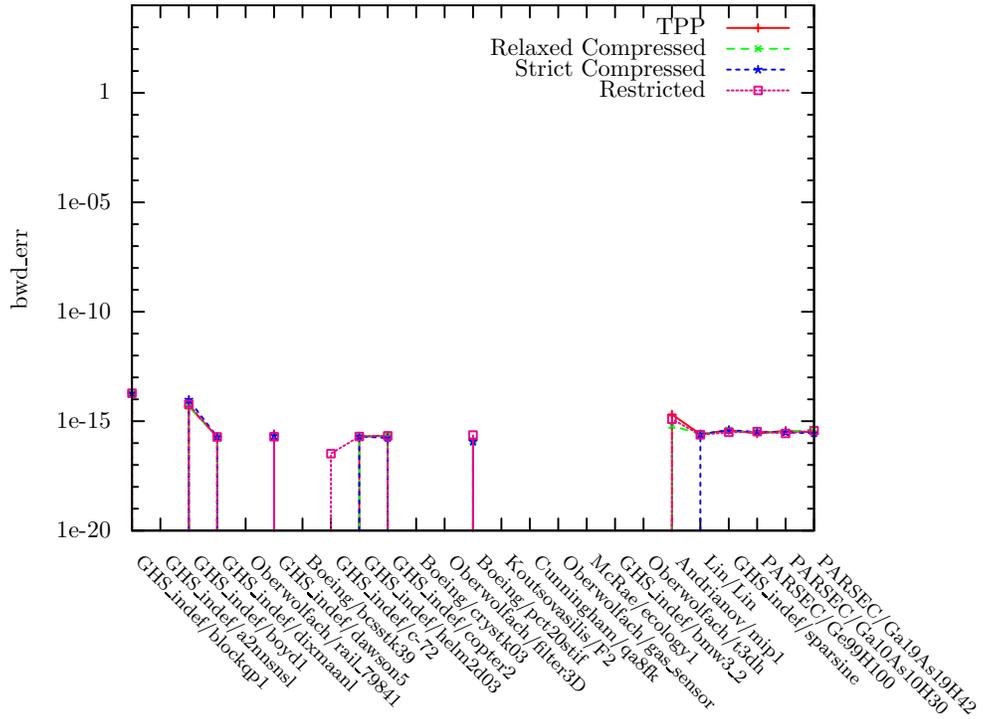

      \caption{ \label{bwd errors}
         Backward errors after iterative refinement (default ordering).
      } 
      \centering
      Test Set 1 \\
      \scalebox{0.90}{\input{general_bwd.tex}} \\
      Test Set 2 \\
      \scalebox{0.90}{\input{delays_bwd.tex}}
   \end{figure}

   For the solution of the indefinite system $\mathcal{A}x=b$, Figure~\ref{bwd errors} plots the
   scaled backward error
   $$
   \mbox{bwd err} = \frac {\|\mathcal{A}x-b\|_{\infty}}
         {\|\mathcal{A}\|_{\infty} \|x\|_{\infty} + \|b\|_{\infty}},
   $$
   after 10 steps of iterative refinement for the pivoting
   strategies described in this paper (used with the default ordering). As expected, they each perform
   adequately on the general problems of Test Set 1. However, on the more
   numerically challenging problems of Test Set 2, the restricted pivoting
   approach fails (that is, iterative refinement fails to converge to a backward
   error smaller than $10^{-14}$),
   while the two numerically stable approaches (threshold
   partial pivoting and strict compressed pivoting) solve every problem to
   machine precision. The relaxed compressed pivoting fails on 3
   problems.
   
\begin{figure}
      \caption{ \label{bwd errors mc80}
         Backward errors after iterative refinement (matching-based ordering).
      } 
      \centering
      Test Set 2 \\
      \scalebox{0.90}{\input{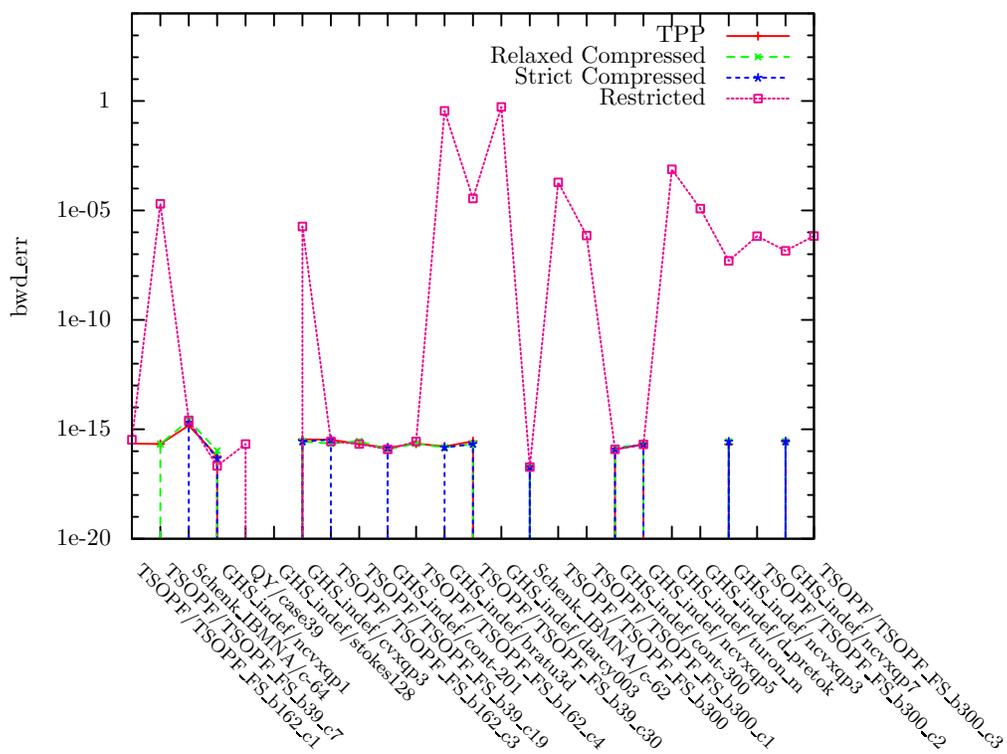}}
   \end{figure}

   By applying a matching-based ordering to the problems in Test Set 2, 
   we obtain the results presented in
   Figure~\ref{bwd errors mc80}. Restricted pivoting still fails to converge to
   an accurate answer for 13 of the 25 problems, but all failures for the relaxed compressed pivoting
   are eliminated. These results demonstrate that, although not backward stable,
   if combined with a matching-based ordering and scaling, relaxed compressed pivoting
   is stable in practice.
   
   Figures~\ref{delay plots} and \ref{delay plots mc80}
   present the numbers of delayed pivots for the default and matching-based
   orderings, respectively. Strict compressed pivoting generally
   results in more delayed pivots than TPP, while for relaxed compressed pivoting the
   number is the same or fewer than for TPP.
   We remark that a small number of delayed pivots
   (typically less than 1000  for problems of
   the size used in our tests) has no significant effect
   on performance.
   The problems in Test Set 1 with the default ordering
   demonstrate the weakness of strict compressed pivoting.
   With the exception of problems GHS\_indef/c-72 and GHS\_indef/bmw3\_2, threshold
   partial pivoting and relaxed compressed pivoting give few (if any) delayed pviots.
   But the stricter pivot selection of strict compressed pivoting results in
   the generation of over 1000 times more delayed pivots
   for some problems. In performance terms, 
   for 5 of the 25 problems in Test Set 1 the {\tt HSL\_MA97} time 
   using strict compressed pivoting
   is more than twice that of using threshold partial
   pivoting (sometimes more than four times greater).
   Of course, for numerically straightforward problems such as these, the
   number of delayed pivots can be reduced by using a smaller threshold
   parameter $u$ without compromising stability and this is an option
   that may want to be considered with strict compressed pivoting. Using a matching-based ordering also
   substantially reduces the number of delayed pivots (but may involve
   more operations and greater fill-in).

   \begin{figure}
      \caption{ \label{delay plots}
         The number of additional delays generated by compressed pivoting compared
         with the number generated by threshold partial pivoting (default ordering).
      } 
      \centering
      Test Set 1 \\
      \scalebox{1.20}{\input{general_ndelay.tex}} \\
      Test Set 2 \\
      \scalebox{1.20}{\input{delays_ndelay.tex}}
   \end{figure}
   
   \begin{figure}
      \caption{ \label{delay plots mc80}
         The number of additional delays generated by compressed pivoting compared
         with the number generated by threshold partial pivoting (matching ordering).
      } 
      \centering
      Test Set 1 \\
      \scalebox{1.20}{\input{general_ndelay_mc80.tex}} \\
      Test Set 2 \\
      \scalebox{1.20}{\input{delays_ndelay_mc80.tex}}
   \end{figure}

   \newpage
   \section{Conclusions} \label{Sec:conclusions}

   Motivated by the need to devise algorithms that communicate as little
   as possible, even if they do slightly more arithmetic operations, 
   we have presented two variants of a new pivoting algorithm 
   for use within a sparse symmetric indefinite direct solver. Our proposed
   variants construct a compressed matrix using a tree reduction algorithm.
   We have shown that this results in 
   better communication properties both practically and asymptotically than
   threshold partial pivoting, without compromising numerical robustness.
   Numerical tests demonstrate over a two times speedup for large problems.
   The strict compressed pivoting algorithm is
   provably stable but at the cost of potentially generating significantly more
   delayed pivots than  threshold partial pivoting. An alternative
   relaxed compressed pivoting  algorithm  avoids this problem, but may not be stable
   on pathological examples. Nonetheless, it is shown that, in combination with
   appropriate scaling and ordering algorithms, it is stable in practice on even
   the most difficult of practical problems.   

   We note that many problems, if well scaled, do not require numerical pivoting, and in such cases the
   try-it-and-see approach suggested by Kim and Eijkhout~\cite{kiei:2012}
   may be more appropriate. However
   for problems where unacceptable growth in the factor entries is detected, our new
   approach offers a
   fast alternative to identifying the minimal set of pivots that must be
   delayed during the factorization. Our future work is to develop software
   that uses such a technique and is targeted at  manycore architectures such as GPUs or Intel's Xeon Phi.

   \section*{Acknowledgements}
   We would like to thank Dianne O'Leary and Iain Duff for their comments on a draft
   of this paper.

   \bibliographystyle{siam}
   \bibliography{ref}

\begin{thebibliography}{10}

\bibitem{adek:2001}
{\sc P.~Amestoy, I.~Duff, J.-Y. L'Excellent, and J.~Koster}, {\em A fully
  asynchronous multifrontal solver using distributed dynamic scheduling}, SIAM
  J. Matrix Analysis and Applications, 23 (2001), pp.~15--41.

\bibitem{agl:1998}
{\sc C.~Ashcraft, R.~Grimes, and J.~Lewis}, {\em Accurate symmetric indefinite
  linear equation solvers}, SIAM J. Matrix Analysis and Applications, 20
  (1998), pp.~513--561.

\bibitem{basd:60}
{\sc D.~Barron and H.~Swinnerton-Dyer}, {\em Solution of simultaneous linear
  equations using magnetic-tape store}, The Computer J., 3 (1960), pp.~28--33.

\bibitem{bebd:2011}
{\sc D.~Becker, M.~Baboulin, and J.~Dongarra}, {\em Reducing the amount of
  pivoting in symmetric indefinite systems}, {T}echnical {R}eport ICL-UT-11-06,
  University of Tennessee, 2011.
\newblock Also INRIA Research Report 7621.

\bibitem{buka:77}
{\sc J.~Bunch and L.~Kaufmann}, {\em Some stable methods for calculating
  inertia and solving symmetric linear systmes}, Mathematics of Computation, 31
  (1977), pp.~163--179.

\bibitem{Davis2011}
{\sc T.~Davis and Y.~Hu}, {\em The {U}niversity of {F}lorida {S}parse {M}atrix
  {C}ollection}, ACM Transactions on Mathematical Software, 38 (2011).
\newblock Article 1, 25 pages.

\bibitem{duer:86}
{\sc I.~Duff, A.~Erisman, and J.~Reid}, {\em {D}irect {M}ethods for {S}parse
  {M}atrices}, Oxford University Press, 1986.

\bibitem{duko:2001}
{\sc I.~Duff and J.~Koster}, {\em On algorithms for permuting large entries to
  the diagonal of a sparse matrix}, SIAM J. Matrix Analysis and Applications,
  22 (2001), pp.~973--996.

\bibitem{dupr:2007}
{\sc I.~Duff and S.~Pralet}, {\em Towards a stable static pivoting strategy for
  the sequential and parallel solution of sparse symmetric indefinite systems},
  SIAM J. Matrix Analysis and Applications, 29 (2007), pp.~1007--1024.

\bibitem{grdx:2011}
{\sc L.~Grigori, J.~Demmel, and H.~Xiang}, {\em {CALU}: A communication optimal
  {LU} factorization algorithm}, SIAM J. Matrix Analysis and Applications, 32
  (2011), pp.~1317--1350.
\newblock Also LAPACK Working Note 266.

\bibitem{gujk:2001}
{\sc A.~Gupta, M.~Joshi, and V.~Kumar}, {\em {WSMP}: {A} high-performance
  serial and parallel sparse linear solver}, Technical Report RC 22038 (98932),
  IBM T.J. Watson Research Center, 2001.
\newblock \url{http://www.cs.umn.edu/~agupta/doc/wssmp-paper.ps}.

\bibitem{hosc:2010c}
{\sc J.~Hogg and J.~Scott}, {\em An indefinite sparse direct solver for large
  problems on multicore machines}, Technical Report RAL-TR-2010-011, Rutherford
  Appleton Laboratory, 2010.

\bibitem{hosc:2011a}
\leavevmode\vrule height 2pt depth -1.6pt width 23pt, {\em {HSL\_MA97}: a
  bit-compatible multifrontal code for sparse symmetric systems}, {T}echnical
  {R}eport RAL-TR-2011-024, Rutherford Appleton Laboratory, 2011.

\bibitem{hosc:2013}
\leavevmode\vrule height 2pt depth -1.6pt width 23pt, {\em A study of pivoting
  strategies for tough sparse indefinite systems}, ACM Transactions on
  Mathematical Software, to appear (2013).
\newblock See also Technical Report RAL-TR-2012-009, Rutherford Appleton
  Laboratory.

\bibitem{hsl:2013}
{\sc HSL}, {\em A collection of {F}ortran codes for large-scale scientific
  computation}, 2013.
\newblock \url{http://www.hsl.rl.ac.uk/}.

\bibitem{kiei:2012}
{\sc K.~Kim and V.~Eijkhout}, {\em A parallel sparse direct solver via
  hierarchical dag scheduling}, {T}echnical {R}eport TR-12-05, Texas Advanced
  Computing Centre, 2012.

\bibitem{lide:98}
{\sc X.~Li and J.~Demmel}, {\em Making sparse {Gaussian} elimination scalable
  by static pivoting}, in Proceedings of the 1998 ACM/IEEE conference on
  Supercomputing, IEEE Computer Society, 1998, pp.~1--17.

\bibitem{Reid2008}
{\sc J.~Reid and J.~Scott}, {\em An efficient out-of-core sparse symmetric
  indefinite direct solver}, {T}echnical {R}eport RAL-TR-2008-024, Rutherford
  Appleton Laboratory, 2008.

\bibitem{resc:2011}
\leavevmode\vrule height 2pt depth -1.6pt width 23pt, {\em Partial
  factorization of a dense symmetric indefinite matrix}, ACM Transactions on
  Mathematical Software, 38 (2011).

\bibitem{scga:2006}
{\sc O.~Schenk and K.~G\"artner}, {\em On fast factorization pivoting methods
  for symmetric indefinite systems}, Electronic Transactions on Numerical
  Analysis, 23 (2006), pp.~158--179.

\bibitem{scwh:2007}
{\sc O.~Schenk, A.~W\"achter, and M.~Hagemann}, {\em Matching-based
  preprocessing algorithms to the solution of saddle-point problems in
  large-scale nonconvex interior-point optimization}, Computer Optimization and
  Applications, 36 (2007), pp.~321--341.

\bibitem{sore:85}
{\sc D.~Sorenson}, {\em Analysis of pairwise pivoting in {G}aussian
  elimination}, IEEE Transactions on Computers, C-34 (1985), pp.~274--278.

\bibitem{trsc:90}
{\sc L.~Trefethhen and R.~Schreiber}, {\em Average-case stability of {G}aussian
  elimination}, SIAM J. Matrix Analysis and Applications, 11 (1990),
  pp.~335--360.

\end{thebibliography}

   \newpage
   \appendix

   \section{Communication analysis} \label{appendix:comm anal}
   \subsection{General results}
   \subsubsection{Serial factorization}
   The serial factorization of a $n\times p$ supernodal matrix forms the basis of much of
   our analysis of the parallel algorithms. To simplify our analysis, we will
   assume henceforth that all $2\times 2$ pivots are accepted immediately and
   hence there is no need to apply permutations.

   For each of $p/2$ pivots, the serial factorization performs the following
   steps. Counts are given for the $i$-th pivot.
   \begin{itemize}
      \item Find column maxima. 2 columns each with $(n-2i)$ entries below the
         pivot require $2(n-2i-1)$ operations.
      \item Test pivot for acceptability. $18$ operations.
      \item Apply (inverse of) pivot to 2 columns. $4$ operations per row for
         total of $4(n-2i)$ operations.
      \item Update the $(p-2i)$ columns to the right of pivotal columns.
         Because of trapezoidal nature of the supernodal contribution to the
         factors, the number of entries to update is
         $$
            (n-p)(p-2i) + \sum_{j=1}^{p-2i} j = \frac{1}{2}(p-2i)(2n-p-2i+1).
         $$
         Each entry updated requires two fused multiply-add operations (one for
         each column of the pivot).
   \end{itemize}
   Summing across all pivots, we obtain the following operation count.
   \begin{eqnarray*}
   \mathrm{TPP}_{\mathrm{ops}}(n,p) & = & \sum_{i=1}^{p/2}\left[ 2(n-2i-1) + 18 + 4(n-2i) + (2n-p-2i+1)\right] \\
      & = & (16+p-p^2+6n+2np)\frac{p}{2} + (-14-4n)\sum_{i=1}^{p/2}i + 4\sum_{i=1}^{p/2}i^2 \\
      & = & (16+p-p^2+6n+2np)\frac{p}{2} + \frac{1}{2}(-14-4n)\frac{p}{2}(\frac{p}{2}+1) + 4\frac{1}{6}\frac{p}{2}(\frac{p}{2}+1)(p+1) \\
      & = & \frac{29}{6}p -\frac{3}{4}p^2 -\frac{1}{3}p^3 + 2np + \frac{1}{2}np^2.
   \end{eqnarray*}

   \subsubsection{Reduction on a tree}
   To perform the communication analysis it is  necessary to have a model
   for reduction. Consider performing simultaneous reduction of $k$ values on
   a binary tree, where only a single processor needs the final result.
   A binary tree across $P$ processors has $1+\log_2 P$
   levels and therefore requires $\log_2 P$ messages to be sent. At each
   non-leaf node of the tree, $k$ comparison operations performed.
   Given that there are $P$ leaf nodes, the number of non-leaf nodes is
   $$
   \sum_{i=1}^{\log_2 P} P2^{-i} = (1-2^{-\log_2 P})P = (P-1).
   $$
   Each non-leaf node has
   $2k$ words of information communicated to it. Hence,
   \begin{eqnarray*}
      \mathrm{Red}_{\mathrm{ops}}(k) & = & (P-1)k \\
      \mathrm{Red}_{\mathrm{msg}}(k) & = & \log_2 P \\
      \mathrm{Red}_{\mathrm{bw}}(k) & = & 2(P-1)k.
   \end{eqnarray*}

   \subsection{Threshold Partial Pivoting}
   We proceed to calculate theoretical bounds on the communication and
   computation for Variants A and B described in Section~\ref{TPP:variants}.

   To simplify the analysis for Variant A, 
   we assume that the first $p$ rows ($A_{11}$)
   all reside on a single processor. Recall that for Variant B, all processors
   are sent a copy of these rows at the start
   (1 message using \\ $\frac{1}{2}(P-1)p(p+1)$ words).

   The significant differences from the serial factorization are in finding
   the column maxima (both variants) and communicating the chosen pivot
   (Variant A) or updating the local copy of $A_{11}$ (Variant B).

   In finding the column maxima for a single pivot, $2(P-1)$
   comparisons are replaced by a global reduction of two values
   (one for each column). This generates no extra operations, but does
   generate an extra $\log_2 P$ messages containing a total of $4(P-1)$ words
   for each pivot. 

   For Variant A, the pivot and first $(p-2i)$ rows of the pivot columns
   must be communicated to other processors. This requires no additional
   operations, but requires 1 message and $2P(p-2i)+3$ words for pivot $i$
   ($1\le i\le p/2$). Hence,
   \begin{eqnarray*}
      \mathrm{TPP}_{\mathrm{ops}}^A(n,p) &=& \mathrm{TPP}_{\mathrm{ops}}(n,p) \\
         %&=& \frac{29}{6}p -\frac{3}{4}p^2 -\frac{1}{3}p^3 + 2np + \frac{1}{2}np^2\\
         &=& O(p^3+np^2) \\
      \mathrm{TPP}_{\mathrm{msgs}}^A(n,p) &=&  p + \frac{1}{2}p\log_2 P \\
         &=& O(p\log P) \\
      \mathrm{TPP}_{\mathrm{bw}}^A(n,p) &=&
         \sum_{i=1}^{p/2}\left[ 4(P-1) + 2P(p-2i)+3 \right] \\
         &=& (-\frac{1}{2} + 2P + Pp)p - 2P(p/2)(p/2+1) \\
         &=& -\frac{1}{2}p + \frac{1}{2}Pp(p+2) \\
         &=& O(Pp^2).
   \end{eqnarray*}

   For Variant B, the leading $p\times p$ submatrix must be updated on every
   processor, incurring an additional $(P-1)\mathrm{TPP}_{\mathrm{ops}}(p,p)$
   operations. This gives
   \begin{eqnarray*}
      \mathrm{TPP}_{\mathrm{ops}}^B(n,p) & = & \mathrm{TPP}_{\mathrm{ops}}(n,p) + (P-1)\mathrm{TPP}_{\mathrm{ops}}(p,p)\\
         & = & \mathrm{TPP}_{\mathrm{ops}}(n,p) + (P-1)(\frac{29}{6}p + \frac{5}{4}p^2 +\frac{1}{6}p^3) \\
         %&=& -2p^2 -\frac{1}{2}p^3 + 2np + \frac{1}{2}np^2 + P(\frac{29}{6}p + \frac{5}{4}p^2 +\frac{1}{6}p^3) \\
         &=& O(p^3+np^2+Pp^3) \\
      \mathrm{TPP}_{\mathrm{msgs}}^B(n,p) & = & 1 + \frac{1}{2}p\log_2 P \\
         &=& O(p\log P) \\
      \mathrm{TPP}_{\mathrm{bw}}^B(n,p) &=&
         \frac{1}{2}(P-1)p(p+1) + \sum_{i=1}^{p/2}\left[ 4(P-1) \right] \\
         &=& \frac{1}{2}(P-1)p(p+1) + 2(P-1)p \\
         &=& -\frac{1}{2}p(3 -p)+ \frac{1}{2}Pp(p+5) \\
         &=& O(Pp^2).
   \end{eqnarray*}

   \subsection{Restricted Pivoting}
   To simplify the analysis, we again assume $A_{11}$ resides on
   a single processor. First,  $A_{11}$ is
   factorized using serial threshold partial pivoting that requires no
   communication and has an operation count
   $$
      \mathrm{TPP}_{\mathrm{ops}}(p,p) = \frac{29}{6}p + \frac{5}{4}p^2 +\frac{1}{6}p^3.
   $$

   The factor $L_{11}$ is then communicated to all other processors, requiring 1 message
   per processor and $\frac{1}{2}(P-1)p(p+1)$ words of bandwidth. Each processor
   then applies this matrix to its own data. The number of operations required
   is
   \begin{eqnarray*}
      (n-p)\sum_{i=1}^{p/2}\left[4+2(p-2i)\right] & = & (n-p)(2p+p^2-2\frac{p}{2}(\frac{p}{2}+1))\\
      & = & \frac{1}{2}(n-p)p(2+p).
   \end{eqnarray*}

   Overall, we thus have:
   \begin{eqnarray*}
      \mathrm{Restrict}_{\mathrm{ops}}(n,p) & = & \frac{29}{6}p + \frac{5}{4}p^2 +\frac{1}{6}p^3 + \frac{1}{2}(n-p)p(2+p) \\
         & = & \frac{29}{6}p +\frac{1}{4}p^2-\frac{1}{3}p^3+np+\frac{1}{2}np^2 \\
         & = & \mathrm{TPP}_{\mathrm{ops}}(n,p) - p(n-p) \\
         & = & O(p^3+np^2) \\
      \mathrm{Restrict}_{\mathrm{msg}}(n,p) &=& 1 \\
      \mathrm{Restrict}_{\mathrm{bw}}(n,p) &=& -\frac{1}{2}p(p+1) + \frac{1}{2}Pp(p+1). \\
   \end{eqnarray*}

   \subsection{Compressed Pivoting}
   Analysis of compressed pivoting (Section~\ref{Sec:STP}) follows that of
   restricted pivoting, except
   for the treatment of the leading submatrix $A_{11}$. The compressed matrix
   $C$ is first assembled. In the strict algorithm, each
   processor scans the rows  assigned to it and places them in the relevant set $J_i$. For each
   of $(n-p)$ rows this requires $p$ absolute value operations and $(p-1)$
   comparisons. Assuming the compressed matrix is initialized to zero,
   filling it requires one additional comparison for each of the $(n-p)p$ matrix
   entries. A reduction on $p^2$ values is then performed. In the relaxed
   algorithm, each processor scans each column of the rows assigned to it to find the largest local entry and
   compares that to the current largest for that row. This process requires a total
   of $(n-p)$ absolute value operations and $(n-p-P)+P$ operations per column.
   A modified reduction is then performed that uses
   $\mathrm{Red}_{\mathrm{ops}}(p)$ operations but
   $\mathrm{Red}_{\mathrm{bw}}(p^2)$ words of bandwidth. We summarise these
   results in the following table. \\
   \begin{table}[h]
      \centering
      \caption{
         Cost for construction of compressed matrix $C$.
      }
      \begin{tabular}{l|c|c|c}
                      & operations & messages & bandwidth \\
         \hline
         Strict   & $(n-p)(3p-1) + (P-1)p^2$ & $\log_2 P$ & $2(P-1)p^2$ \T \\
         Relaxed & $(n-p)(2p) + (P-1)p$ & $\log_2 P$ & $2(P-1)p^2$ \\
      \end{tabular}
   \end{table}

   Having constructed $C$, the modified factorization is  performed. The strict algorithm
   requires additional operations to work with the absolute value update
   operation on the $p\times p$ matrix $C$. This can be done efficiently
   by storing an extra copy of $A_{11}$ with the
   absolute value operation applied. The number of operations required is
   $$
      \mathrm{TPP}_{\mathrm{ops}}(2p,p) + \frac{1}{2}p(p+1)
      = \frac{32}{6}p + \frac{15}{4}p^2  + \frac{2}{3}p^3.
   $$
   As the relaxed algorithm uses an unmodified TPP algorithm, the
   number of operations it requires is
   $$
      \mathrm{TPP}_{\mathrm{ops}}(2p,p) = \frac{29}{6}p + \frac{13}{4}p^2  + \frac{2}{3}p^3.
   $$

   Summing with the application of $L_{11}$ to $A_{21}$, we obtain the
   following counts,
   \begin{eqnarray*}
      \mathrm{Compressed}_{\mathrm{ops}}^{WC}(n,p) & = & (n-p)(3p-1)+(P-1)p^2 + \frac{32}{6}p + \frac{15}{4}p^2 + \frac{2}{3}p^3 + \frac{1}{2}(n-p)p(2+p) \\
         & = & Pp^2 + \frac{19}{3}p - \frac{5}{4}p^2 + \frac{1}{6}p^3 - n + 4np + \frac{1}{2}np^2 \\
         & = & \mathrm{TPP}_{\mathrm{ops}}(n,p) + \frac{1}{2}p((p-1)p+3) + n(2p-1) + Pp^2 \\
         & = & O(Pp^2 + p^3 + np^2) \\
      \mathrm{Compressed}_{\mathrm{ops}}^{AC}(n,p) & = & (n-p)(2p)+(P-1)p + \frac{29}{6}p + \frac{13}{4}p^2 + \frac{2}{3}p^3 + \frac{1}{2}(n-p)p(2+p) \\
         & = & Pp + \frac{23}{6}p + \frac{1}{4}p^2 + \frac{1}{6}p^3 + 3np + \frac{1}{2}np^2 \\
         & = & \mathrm{TPP}_{\mathrm{ops}}(n,p) + \frac{1}{2}p((p+2)p-2) + (n+P)p \\
         & = & O(Pp + p^3 + np^2). \\
   \end{eqnarray*}
   The strict and relaxed algorithms have the same
   communication pattern. Summing the compressed matrix construction with the later distribution of $L_{11}$ yields
   \begin{eqnarray*}
      \mathrm{Compressed}_{\mathrm{msg}}(n,p) & = & 1 + \log_2 P \\
                                              & = & O(\log P) \\
      \mathrm{Compressed}_{\mathrm{bw}}(n,p) & = & 2(P-1)p^2 + \frac{1}{2}(P-1)p(p+1) \\
                                             & = & -\frac{1}{2}p(5p+1) + \frac{1}{2}Pp(5p+1) \\
                                           & = & O(Pp^2).
   \end{eqnarray*}
\end{document}